\documentclass{art_ar}

\usepackage{amsmath,amssymb}

\newcommand{\R}{\mathbb{R}}
\newcommand{\RR}{\mathbb{R}}
\newcommand{\N}{\mathbb{N}}
\newcommand{\EX}{\mathbf{E}}

\newcommand{\X}{\overline{X}^{(n)}}
\newcommand{\XM}{\widehat{X}^{(n)}}

\newcommand{\YM}{\widehat{Y}}

\newcommand{\XL}{{\widetilde{X}}^{(n)}}
\newcommand{\YL}{{\widetilde{Y}}^{(n)}}
\newcommand{\BL}{{\widetilde{B}}^{(n)}}

\newcommand{\HH}{\mathcal{H}}
\newcommand{\EE}{\mathcal{E}}
\newcommand{\sk}{\mathbb{D}}

\newcommand{\fin}{
\begin{flushright}
\mbox{$\square$}
\end{flushright}
\noindent }

\begin{document}
\begin{frontmatter}
\title{    Optimal pointwise approximation of stochastic differential equations driven  by fractional
  Brownian motion             }
\author{Andreas Neuenkirch}
\corauth[cor1]{neuenkirch@math.uni-frankfurt.de}
\address{Johann Wolfgang Goethe-Universit\"at Frankfurt, \\
Institut f\"ur Mathematik,
Robert-Mayer-Str. 10,\\
D-60325 Frankfurt am Main,
Germany 
}

\thanks{The author is supported by the DFG-project ``Pathwise numerical analysis of stochastic evolution equations''.}

\begin{abstract}
We study the  approximation  of stochastic differential equations driven by a fractional
Brownian motion with Hurst parameter $H>1/2$.  For the
mean-square error at a single point  we derive the optimal
rate of convergence that can be achieved by any approximation method using an equidistant discretization of the
driving fractional Brownian motion. 
We find that there are mainly two cases: either the solution can be approximated perfectly or the
best possible rate of convergence is $n^{-H-1/2},$ where $n$ denotes the number of evaluations of the fractional Brownian
motion.
In addition, we present an implementable approximation scheme that obtains the optimal rate of convergence in
the latter case.
\end{abstract}


\end{frontmatter}

\section{Introduction}
In this article, we study the pathwise approximation of the equation 
\begin{align}\label{eq1} dX_t&=a(X_t) \, dt + \sigma(X_t)dB_t,\qquad t\in [0,1], \\  
 X_0&=x_0, \nonumber  \end{align} 
where $x_0 \in \R$ and $(B_t)_{t \in [0,1]}$ is a fractional Brownian motion with Hurst parameter $H \in (1/2,1)$. 
 Equation (\ref{eq1}) is understood as a pathwise Riemann-Stieltjes  integral equation, see e.g.  \cite{lpma}, \cite{nua_saus} and \cite{zaehle}. 
Recent applications of  stochastic differential equations driven by fractional Brownian motion include e.g. the noise simulation in electronic circuits (\cite{denk}), the modelling of the subdiffusion
of a protein molecule (\cite{protein2}) and the pricing of weather derivatives (\cite{benth}, \cite{brody}).

The type of approximation methods we are concerned with produce an approximation to $X_1$ using an equidistant discretization of the driving fractional Brownian motion, that is $$ B_{1/n} \, , \,  B_{2/n} \, , \,  \ldots  \, , \, B_{1}.$$ The error of such an approximation method will be measured by the mean square norm $(\EX | \cdot |^{2} )^{1/2}$. 
Thus, the best possible approximation method of the above type is clearly the conditional expectation 
$$ \X_1 = \EX ( \, X_1  \, | \,  B_{1/n} , B_{2/n}, \ldots, B_{1} \, ).$$
For stochastic differential equations driven by Brownian motion, the optimal pointwise approximation (in the mean-square sense) of the solution
is a well studied problem, also for non-equidistant and adaptive discretizations. See, e.g., \cite{cam_clark}, \cite{new1}, \cite{new2}, \cite{cas_gaines}, \cite{cam_hu} and \cite{MG}.
However, for stochastic differential equations driven by fractional Brownian motion there are only a few known results for 
mean-square approximation, mainly for equations with additive noise (\cite{joc}, \cite{diss}) or with a linear diffusion coefficient (\cite{MS}).
In  \cite{joc} and \cite{diss} the Euler and a Wagner-Platen-type method for equations with additive noise are studied and their exact rates of convergence are given, while
in \cite{MS} the convergence order of an Euler-type method for a quasi-linear Skorohod-type equation is determined. 
Moreover, the asymptotic error distribution of several approximation schemes for equation (\ref{eq1}) is derived in  \cite{NN} and \cite{GN}. 
 
Throughout this article, we will impose rather strong assumptions on the drift- and diffusion coefficients, which has technical reasons. See Remark (i) in Section \ref{remarks_section5}  for a detailed discussion. 
In particular $a$ and $\sigma$ are supposed to be bounded and $\sigma$ is also strictly positive.
Nevertheless, we think that this article will give a lot of structural insight in the approximation of  stochastic differential equations driven by fractional Brownian motion.

If the drift and diffusion coefficient commute, i.e. $a'\sigma- a \sigma' =0$, then  it is easy 
to show that $X_1$ does not depend on the whole process $(B_{t})_{t \in [0,1]}$ but only  on  $B_1$.
This implies that 
\begin{align} (\EX |X_1- \X_1 |^{2} )^{1/2} = 0.\end{align}
More generally, we can show the following relation:  
\begin{center} There exists a mapping $g \in C^{2}([0,1] \times \R;  \R)$ such  that 
 $X_t=g(t,B_t)$  holds \\ for all $t\in [0,1]$ almost surely, if and only if $(a'\sigma- a \sigma')(x) =0$ for all $x \in \RR$. \end{center}
This extends  a particular case of a well known result for  stochastic differential equations driven by Brownian motion. See e.g. \cite{yamato}.
 
Otherwise, if $a$ and $\sigma$ do not commute, we have the following upper and lower error bound for the error of the conditional expectation
\begin{align} \liminf_{n \rightarrow \infty} \,  n^{H+1/2} \,  (\EX |X_1- \X_1|^{2} )^{1/2} 
 \, \geq  \,  \alpha_H \left(\int_{0}^{1} |\EX \mathcal{Y}_t|^{2} \, dt \right)^{1/2} \label{low_bound} \end{align}
and
\begin{align}  \limsup_{n \rightarrow \infty} \, n^{H+1/2}  \, (\EX |X_1- \X_1 |^{2} )^{1/2} 
 \, \leq  \,  \beta_H \left(\int_{0}^{1} \EX| \mathcal{Y}_t|^{2} \, dt \right)^{1/2},  \label{upp_bound}      \end{align}
where $\alpha_H$ and $\beta_H$ are two constants, depending only on $H$ and the random weight function 
$(\mathcal{Y}_t)_{t \in [0,1)}$ is  given by
$$ \mathcal{Y}_t = (a\sigma' - a'\sigma)(X_t) \exp \left( \int_{t}^{1} a'(X_s) \,ds 
+ \int_{t}^{1} \sigma'(X_s) \, dB_s \right), \qquad t \in [0,1].$$
If the constant on the right hand side of equation (\ref{low_bound}) does not vanish, then
the conditional expectation has exact rate of convergence $n^{-H-1/2}$. This is satisfied  for example, 
if $(a\sigma' - a' \sigma)(x_0) \neq 0$.  Consequently, in this case
there is no approximation method using an equidistant discretization of the driving fractional
Brownian motion that can obtain a better rate of convergence than $n^{-H-1/2}$.

The conditional expectation is clearly in general not an implementable method for the approximation of stochastic differential equations driven by fractional Brownian motion.
Therefore, we also consider here  an extension of  McShane's method for
Stratonovich SDEs, see \cite{macshane} and \cite{MG} for a related scheme for It\^{o} SDEs.  Our extension of the McShane method is  defined by 
$\XM_{0}=x_{0}$ and \begin{align*} \XM_{k+1}=\XM_{k} &+ a(\XM_{k}) \Delta + \sigma(\XM_{k})\Delta_k B + \frac{1}{2} \sigma \sigma'(\XM_k) (\Delta_k B)^{2} \\ &+\frac{1}{2}(a\sigma' + a'\sigma)(\XM_k) \Delta_k B \Delta + \frac{1}{2} aa'(\XM_k) \Delta^{2} \\ & + \frac{1}{6} (\sigma^{2} \sigma'' +  \sigma(\sigma')^{2})(\XM_k)(\Delta_k B)^{3}  \end{align*} for $k=0, \ldots, n-1$, where $\Delta=1/n$ and  $\Delta_k B = B_{(k+1)\Delta} - B_{k \Delta}$.

We show that
\begin{align} \lim_{n \rightarrow \infty} \, n^{H+1/2}  \, (\EX |X_1- \XM_n |^{2} )^{1/2} 
\,=  \,  \beta_H \left(\int_{0}^{1} \EX| \mathcal{Y}_t|^{2} \, dt \right)^{1/2}. \end{align}
Hence this implementable approximation method has
exact rate of convergence  $n^{-H-1/2}$, if $a$ and
$\sigma$ do not commute, and thus obtains the same convergence rate as the conditional expectation  in this case.

The article is structured as follows. In the next section we recall some facts about fractional Brownian
motion and stochastic differential equations driven by fractional Brownian motion.
In the Sections 3 and 4, we state and prove our results for the error of
the conditional expectation, while McShane's method is considered in
section 5. A technical proof of an auxiliary result is postponed to the Appendix.

\section{Preliminaries}
\subsection{Fractional Brownian motion}\label{pre_1}
Let $B=(B_t)_{t \in  [0,1]}$ be a fractional Brownian motion with Hurst parameter $H \in (0,1)$ defined on a complete probability space $(\Omega, \mathcal{F}, P)$, i.e., $B$ is  a continuous centered Gaussian process with covariance function 
$$ R_{H}(s,t)=\frac{1}{2}(s^{2H}+t^{2H}-|t-s|^{2H}),\quad s,t\in [0,1]. $$ 
For $H=1/2$, $B$ is a standard Brownian motion, while for $H\neq 1/2$, it is neither a semimartingale nor a Markov process. Moreover, it holds 
$$ ( \EX |B_t - B_s|^2)^{1/2}=|t-s|^H,\quad s,t\in [0,1], $$ 
and almost all sample paths of $B$ are H\"older continuous of any order $\lambda \in (0,H)$. 

Let us give a few facts about the Gaussian structure of fBm for $H>1/2$ and its Malliavin derivative process, following \cite{nual2} and Chapters 1.2 and 5.2  in \cite{nual1}.
Let $\EE$ be the set of step functions on $[0,1]$ and consider the Hilbert space $\HH$ defined as the closure of $\EE$ with respect to the scalar product 
$$\langle 1_{[0,t]}, 1_{[0,s]}\rangle_{\HH}\; =\; R_H(t,s), \quad s,t \in [0,1].$$
The mapping $1_{[0,t]} \mapsto B_t$ can be extended to an isometry between   $\HH$ and its associated Gaussian space. This  isometry will be denoted by
$\varphi \mapsto B(\varphi)$.
Note that  
\begin{equation}
\label{HH} {\langle \varphi, \rho \rangle}_{\HH}\; = \;\int_0^1  (K_H^*\varphi)(s)(K_H^*\rho)(s)\, ds, \quad \varphi, \rho \in \EE,
\end{equation}
where
\begin{align*}  (K_H^*\varphi)(s)\; =\;  \int_s^1 \varphi(r) \frac{\partial K_H}{\partial r}(r,s)\, dr, \quad \varphi \in \EE, \,\, s \in [0,1],  \end{align*}
with 
$$K_{H}(t,s)=c_{H} s^{1/2-H} \int_{s}^{t} (u-s)^{H-3/2}u^{H-1/2} \, du \, 1_{[0,t)}(s)$$
and
$$c_{H} = \left( \frac{H(2H-1)}{\beta(2-2H,H-1/2)}\right)^{1/2}.$$
Here $\beta$ denotes the Beta function.
Moreover, we have  $L^{1/H}([0,1]) \subset \HH$ and in particular
\begin{align} \langle \varphi, \psi \rangle_{\HH} =  \alpha_{H} \int_{0}^{1} \int_{0}^{1} \varphi(r) \psi(u) |r-u|^{2H-2} \, dr \, du \label{iso_simple} \end{align}
 for $\varphi, \psi \in L^{1/H}([0,1]) $ with $\alpha_H = H(2H-1)$. (For a  characterization of the space $\HH$ in terms of distributions see \cite{jolis}.)

\vspace{2mm}

Let $f : \RR^n \rightarrow \RR$ be a smooth function with compact support and consider the random variable 
$F = f(B_{t_{1}}, \ldots, B_{t_n})$  with $t_{i} \in [0,1]$  for $i=1, \ldots, n$. The derivative process of $F$ is the 
element of $L^2(\Omega; \HH)$ defined by
$$D_sF\; =\; \sum_{i =1}^n \frac{\partial f}{\partial x_i}(B_{t_{1}}, \ldots, B_{t_n})
1_{[0,t_i]}(s), \quad s \in [0,1].$$
In particular $D_s B_t = 1_{[0,t]}(s)$. As usual, $\sk^{1,2}$ is 
the closure of the set of smooth random variables with respect to the norm 
$$\| F\|_{1,2}^2 \, = \,  \EX |F|^2  + \EX \| DF\|_{\HH}^2 .$$

If $F_{1},F_{2} \in \sk^{1,2}$ such that $F_{1}$ and $\|D F_{1}\|_{\HH}$ are bounded, then $F_{1}F_{2} \in \sk^{1,2}$ and 
\begin{align}D F_{1}F_{2} = F_{2} DF_{1} + F_{1} D F_{2}  \label{prod_rule}\end{align} Moreover, recall also the following  chain rule: For $F \in \sk^{1,2}$ and $g \in C^{1}(\RR)$ with bounded derivative we have $g(F) \in \sk^{1,2}$ and
\begin{align}\label{chain_rule_md}
D g(F)= g'(F) \, DF.
\end{align}

The divergence operator $\delta$ is the adjoint of the derivative operator. If a random variable $u \in L^{2}(\Omega; \HH)$
belongs to the domain  $\operatorname{dom}(\delta)$ of the divergence operator denoted, then  $\delta(u)$ is defined by the duality relationship
\begin{align*}  \EX (F \delta(u))= \EX \langle D F, u \rangle_{\HH} \end{align*} for every $F \in \sk^{1,2}$. \\

If $u=(u_t)_{t \in [0,1]}$ is a stochastic process with H\"older continuous sample paths of order $\lambda > 1-H$, 
then the Riemann-Stieltjes integral with respect to $B$ is well defined. If $u$
 moreover satisfies 
$u_t \in \sk^{1,2}$ for all $t \in [0,1]$ and
\begin{align*}\sup_{s \in [0,1]} \EX |u_s|^{2} + \sup_{r,s \in [0,1]} \EX |D_{r}u_s|^{2} < \infty,  \end{align*}
then we have $u \in \textrm{dom}(\delta)$ and the relation
\begin{align}\label{rel} \int_{0}^{1} u_t \, dB_t = \delta(u) + \alpha_H \int_{0}^{1} \int_{0}^{1} D_s u_t |s-t|^{2H-2} \, ds \, dt \end{align}
holds. 
For the Skorohod integral of the process $u$ we have the isometry 
\begin{align}
\label{iso_work}\EX \left|  \delta (u) \right|^{2} &= \, \alpha_{H} \int_{[0,1]^{2}} \EX \,u_s u_r |s-r|^{2H-2} \, dr \, ds \\ &
 \,\,  + \alpha_{H}^{2} \int_{[0,1]^{4}}  \EX D_{r}u_s D_{r'}u_{s'} |r-s'|^{2H-2} |r'-s|^{2H-2} dr dr' ds ds'
\nonumber.\end{align}

\vspace{2mm}

In what follows, we will also the require the Wiener-Chaos decomposition of  a random variable $F \in \sk^{1,2}$: 
Let $H_{n}$ the $n$-th Hermite polynomial, $n \in \N$, and denote by $\mathcal{C}_{n}$ the closed linear subspace of $L^{2}(\Omega, \mathcal{F}, P)$ generated by the random variables $ \{H_{n}(B(\varphi)), \varphi \in \HH, \| \varphi \|_{\HH}=1 \}$. $\mathcal{C}_{0}$ will be the set of constants. Furthermore, denote by $\mathcal{G}$ the $\sigma$-algebra generated by the random variables $B(\varphi), \varphi \in \mathcal{H}$. Then the space $L^{2}(\Omega, \mathcal{G}, P)$ can be decomposed into the infinite orthogonal sum of the subspaces $\mathcal{C}_{n}$:
$$ L^{2}(\Omega, \mathcal{G}, P) = \oplus_{n=0}^{\infty}\,  \mathcal{C}_{n}.$$
Moreover, denote by $\mathcal{J}_{1}: L^{2}(\Omega) \rightarrow \mathcal{C}_{1}$ the projection to the first chaos.
If $F \in \sk^{1,2}$ such that $\sup_{s \in [0,1]} \EX |D_s F|^{2} < \infty$, then we have
\begin{align} \label{chaos} \mathcal{J}_1(F)= \delta(\EX \, D F) \end{align}
almost surely, which is a straightforward consequence of the transfer principle for fractional Brownian
motion, see e.g. Chapter 5.2 in \cite{nual1} and Stroock's formula, see, e.g., Chapter 6 in \cite{p_mall}.

\subsection{Stochastic differential equations driven by fractional Brownian motion}\label{pre_2}

Throughout this article, we will impose the following assumptions:

\begin{itemize}
\item[(A1)] $H > 1/2$
\item[(A2)] $a\in C^{2}(\RR)$, $\sigma \in C^{3}(\RR)$  with bounded derivatives,
\item[(A3)] $a$, $\sigma$ bounded, $\inf_{x \in \RR} \sigma(x)>0$.
\end{itemize}

The assumption (A3)  is  - a priori - required only for technical reasons. See Section \ref{remarks_section5} for a discussion.

Under the assumptions (A1) and (A2) it is well known that
$$ X_t =x_0 + \int_{0}^{t} a(X_{\tau}) \, d \tau+ \int_{0}^{t} \sigma(X_\tau) \, dB_\tau, \qquad t \in [0,1],$$
which is the integral equation  corresponding to equation (\ref{eq1}), has a unique solution
$X=(X_t)_{t \in [0,1]}$
with
$$ \EX \sup_{t \in [0,1]} |X_t|^{p} < \infty$$
for all $p\geq1$. See e.g.~ \cite{HuNu}, \cite{lpma}. Here, the integral with respect to fractional Brownian 
motion is defined as a Riemann-Stieltjes integral.

Moreover, we have $X_t \in \sk^{1,2}$ for all $t \in [0,1]$ and
\begin{align}
\label{mall_deriv} D_uX_t = \sigma(X_u) \exp \left( \int_u^{t} a'(X_{\tau}) \, d \tau + \int_u^{t} \sigma'(X_{\tau}) \, dB_{\tau} \right) 1_{[0,t]}(u), \, \,\, u,t \in [0,1].
\end{align}
See \cite{nosi}, \cite{nua_saus}.
 
Since the diffusion coefficient is strictly positive due to the assumption (A3) we can use
the Lamperti transformation, which will be an important tool throughout this article.
Define
\begin{align*}   \vartheta(x)= \int_{0}^{x} \frac{1}{\sigma(\xi)} \, d \xi, \qquad x \in \RR. \end{align*} Then $\vartheta: \RR \rightarrow \RR$ is well defined, since $\sigma$ is strictly positive. Moreover, we have
$$ \vartheta'(x) = \frac{1}{\sigma(x)}, \qquad x \in \RR. $$
Note that $\vartheta: \RR \rightarrow \RR$ is strictly monotone, thus the inverse function $\vartheta: \RR \rightarrow \RR$ exists 
and satisfies 
$$ (\vartheta^{-1})'(x)= \sigma(\vartheta^{-1}(x)), \qquad x \in \RR.$$
A straightforward application of the change of variable formula for Riemann-Stieltjes integrals,
 see e.g. \cite{zaehle}, yields that 
\begin{align} Y_t= \vartheta(X_t), \qquad t \in [0,1], \end{align}
is the unique solution of  the stochastic differential equation
\begin{align}\label{red_sde}
dY_t&= g(Y_t) \, dt + dB_t, \qquad t \in [0,1], \\
Y_0 &=  \vartheta(x_{0}), \nonumber
\end{align}
with
\begin{align*} g(x)=\frac{a(\vartheta^{-1}(x))}{\sigma(\vartheta^{-1}(x))}, \qquad x \in \RR.   \end{align*}
Clearly, we also have 
\begin{eqnarray} X_t=\vartheta^{-1}(Y_t), \qquad t \in [0,1]. \label{rep_sol}  \end{eqnarray}

Note that  the mapping  $g: \R \rightarrow \R$ is twice continuously differentiable with bounded derivatives.

Using  (\ref{rep_sol})   we can also give a different representation
of the Malliavin derivative of $X_t$, $t \in [0,1]$, which will be more appropriate for our purposes.  See \cite{rindis} for a similar representation in the case of It\^{o} stochastic differential equations driven by Brownian motion.

\medskip

\begin{prop}\label{rep_mall} We have
\begin{align}
 D_{s}X_t=  \sigma(X_t) \exp \left( \int_{s}^{t}  \left( a'  - \frac{a \sigma'}{\sigma} \right) (X_{\tau}) \, d \tau \right)1_{[0,t]}(s) , \qquad s,t \in [0,1].   \label{rep_mall_formula}
\end{align}
\end{prop}
 
\medskip

{\noindent {\bf Proof.}} Using the chain rule (\ref{chain_rule_md}) we have  
$$ D_{s}X_t = (\vartheta^{-1})'(Y_t) D_{s}Y_t = \sigma(X_t) D_{s}Y_t, \qquad s, t \in [0,1]. $$ 
Since $(Y_t)_{t \in [0,1]}$ is the solution of equation (\ref{red_sde}), applying (\ref{mall_deriv}) yields that 
$$ D_s Y_t= \exp \left( \int_{s}^{t}g'(Y_{\tau}) \, d \tau \right)1_{[0,t]}(s),  \qquad s, t \in [0,1]. $$ 
Moreover, it holds
$$   g'(x)=  a'(\vartheta^{-1}(x)) -\frac{a\sigma'}{\sigma}(\vartheta^{-1}(x)). $$
Thus, the  assertion now follows by (\ref{rep_sol}). 
\fin

In particular, the representation (\ref{rep_mall_formula}) implies that  $(D_{s}X_t)_{s,t \in [0,1]}$ is  a bounded stochastic field.

\section{The degenerated case}
In this section we study under which conditions  the solution of equation (\ref{eq1}) is ``degenerated'' in the following sense: The solution $X_t$ at time $t \in [0,1]$ does not depend on the whole sample path $(B_s)_{s \in [0,t]}$ of the driving fractional Brownian motion up to time $t$, but only on $B_t$.  As  for Stratonovich stochastic differential equations driven by Brownian motion, see \cite{yamato},  this property can be completely characterized in terms of the drift- and diffusion coefficient.

\medskip

\begin{thm}\label{main_1}
There exists  a  mapping $f \in C^{2}([0,1] \times \R ; \R)$ such
that
$$ (X_t)_{t \in [0,1]} =\left( f(t,B_t) \right)_{t \in [0,1]} $$
almost surely, if and only if $(a'\sigma- a \sigma')(x) =0$ for all $x \in \R$. 
\end{thm} 

\medskip

{\noindent  {\bf Proof.}} (i) Suppose that $a'\sigma - a \sigma' =0$. 
Using the Lamperti transformation,  we have 
$X_t= \vartheta^{-1}(Y_t)$, $t \in [0,1]$, with
\begin{align*}
dY_t&= g(Y_t) \, dt + dB_t, \qquad t \in [0,1], \\
Y_0 &=  \vartheta(x_{0})
\end{align*}
and
$$g(x)=\frac{a(\vartheta^{-1}(x))}{\sigma(\vartheta^{-1}(x))}, \qquad x \in \RR.$$
Since 
$$ g'(x)=    a'(\vartheta^{-1}(x)) -\frac{a\sigma'}{\sigma}(\vartheta^{-1}(x)) = \frac{ a'\sigma-a\sigma'}{\sigma}(\vartheta^{-1}(x)), $$
the assumption  $a'\sigma - a \sigma' =0$ implies that
 $$g'(x)=0, \qquad x \in \R.$$ 
Consequently, we have
$$ Y_t = \frac{a(x_0)}{\sigma(x_0)} t  +B_t, \qquad t \in [0,1].$$
Thus, we obtain the representation
$$ X_t = f(t,B_t), \qquad t \in [0,1], $$
with
$$f(t,x)=\vartheta^{-1}\left(   \frac{a(x_0)}{\sigma(x_0)} t  + x\right).$$
Due to the assumptions on $\sigma$, we also have that $f \in C^{2}([0,1] \times \R; \R)$.

(ii) Now assume that 
\begin{align} \label{degn} X_t = f(t,B_t), \qquad t \in [0,1] \end{align}
with $f \in C^{2}([0,1] \times \R; \R)$.
From  (\ref{degn}) and the chain rule for Riemann-Stieltjes integrals, see e.g. \cite{zaehle}, we obtain
\begin{align*}
 X_t =x_0 + \int_{0}^{t} f_t(\tau, B_{\tau}) \,  d \tau +    \int_{0}^{t} f_x(\tau, B_{\tau}) \,  d B_{\tau} , \qquad t \in [0,1],
\end{align*}
almost surely. This yields
\begin{align} \label{z1_gl_z2}
\int_{0}^{t}  f_t(\tau, B_{\tau})  - a(X_{\tau}) \, d \tau = -    \int_{0}^{t} f_x(\tau, B_{\tau}) - \sigma(X_{\tau}) \,  d B_{\tau} , \qquad t \in [0,1],
\end{align}
almost surely. Now set
$$ Z_t^{(1)} = \int_{0}^{t}  f_t(\tau, B_{\tau})  - a(X_{\tau}) \, d \tau, \qquad t \in [0,1],$$
and
$$ Z_t^{(2)} =  \int_{0}^{t} f_x(\tau, B_{\tau}) - \sigma(X_{\tau}) \,  d B_{\tau}, \qquad t \in [0,1].$$
Moreover, define the $\alpha$-variation $V^{\alpha}_n(Z)$ with stepsize $1/n$ of a stochastic process $Z=(Z_t)_{t \in [0,1]}$ by
$$ V^{\alpha}_n(Z) = \sum_{i=0}^{n-1} \left | Z_{(i+1)/n} - Z_{i/n} \right|^{\alpha}. $$
Clearly, we have  by (\ref{z1_gl_z2}) that
$$ V^{\alpha}_n(Z^{(1)}) =  V^{\alpha}_n(Z^{(2)}) $$
almost surely for all $\alpha >0$ and $n \in \N$. Since the process $(Z^{(1)}_t)_{t \in [0,1]}$ is pathwise continuously differentiable,
we have that
$$      V^{1/H}_n(Z^{(1)}) \stackrel{Prob}{\longrightarrow} 0   $$
as $n \rightarrow \infty$.
Moreover, Theorem 1 in \cite{woerner} yields that 
$$ V^{1/H}_n(Z^{(2)}) \, \stackrel{Prob}{\longrightarrow} \,  \EX|B_1|^{1/H} \int_{0}^{1} \left|       f_x(\tau, B_{\tau}) - \sigma(X_{\tau})   \right|^{1/H}  d \tau $$
as $n \rightarrow \infty$.
Thus, we have that 
$$  \int_{0}^{1} \left|       f_x(\tau, B_{\tau}) - \sigma(X_{\tau})   \right|^{1/H}  d \tau = 0 $$
almost surely and it follows
$$ f_x(t, B_{t}) =  \sigma(X_{t}),  \qquad t \in [0,1] ,$$ 
almost surely. This yields in turn that
$$ f_t(t, B_{t}) =  a(X_{t}), \qquad t \in [0,1], $$ 
almost surely.
Since $X_t=f(t,B_t)$, we can write the above two equations also as
$$        f_x(t, B_{t}) =  \sigma(f(t,B_t)), \qquad  t \in [0,1].  $$
and 
$$ f_t(t, B_{t}) =  a(f(t,B_t)), \qquad  t \in [0,1].$$
Since $B_t$ is a centered Gaussian random variable with variance $t^{2H}$, we obtain 
$$  f_x(t, x) = \sigma(f(t,x)) , \qquad  f_t(t, x) =  a(f(t,x)), \quad  \qquad t \in (0,1], \, \, x \in \R,$$
which implies that
$$ 0=f_{xt}(t,x)-f_{tx}(t,x) = a\sigma'(f(t,x))- \sigma a'(f(t,x)), \qquad t \in (0,1], \, \, x \in \R.$$
However, since $\sigma$ does not vanish due to our assumptions, the distribution of  $X_t$ is absolutely continuous with respect to the Lebesgue measure with a strictly positive density for every $t \in (0,1]$. See e.g \cite{nosi}. This implies that the image of the mapping $f$ is the whole real line $\R$. Now, we have finally that 
$$         a\sigma'(x) - \sigma a'(x)=0, \qquad  x \in \R.$$
\fin

Thus,  the solution of equation (\ref{eq1}) is ``degenerated'' if and only if $a$ and $\sigma$ commute in the usual sense of differential geometry.
Since
$$  \X_1 = \EX ( \, X_1  \, | \,  B_{1/n} , B_{1/n}, \ldots, B_{1} \, ),$$
the following corollary is a straightforward consequence  of Theorem 1.

\medskip

\begin{cor}\label{perf_app} Let $(a'\sigma - a \sigma')(x) = 0$ for all $x \in \R$. Then we have
$$ ( \EX |X_1 -\X_1|^{2} )^{1/2}=0$$
for all $n \in \N$.
\end{cor}

\medskip

Hence  the above Corollary implies that $X_1$ can be simulated perfectly - at least theoretically. The mapping $f;[0,1] \times \R \rightarrow  \R $ in Theorem \ref{main_1} is given by 
$$  f_x(t, x) = \sigma(f(t,x)) , \qquad  \qquad f_t(t, x) =  a(f(t,x)), \quad  \qquad t \in (0,1], \, \, x \in \R,$$
with $f(0,x)=x_0$, $x \in \R$.  The solution of this partial differential equation  will be explicitly known only in some  particular cases.

\section{The non-degenerate case.}
In this section, we determine the following  lower  and upper bound for the error of the conditional expectation $\X_1$ in the non-commutative case:

\medskip

\begin{thm} \label{main_2} It holds
\begin{align} \liminf_{n \rightarrow \infty} \,  n^{H+1/2} \,  (\EX |X_1- \X_1|^{2} )^{1/2} 
 \, \geq  \,  \alpha_H \left(\int_{0}^{1} |\EX \mathcal{Y}_t|^{2} \, dt \right)^{1/2} \label{low_bound_formula} \end{align}
and
\begin{align} \limsup_{n \rightarrow \infty} \,  n^{H+1/2} \,  (\EX |X_1- \X_1|^{2} )^{1/2} 
 \, \leq  \,  \beta_H \left(\int_{0}^{1} \EX | \mathcal{Y}_t|^{2} \, dt \right)^{1/2} \label{up_bound_formula} \end{align}
where
$$ \mathcal{Y}_t = (a\sigma' - a'\sigma)(X_t) \exp \left( \int_{t}^{1} a'(X_s) \,ds 
+ \int_{t}^{1} \sigma'(X_s) \, dB_s \right), \qquad t \in [0,1]$$ and 
$\alpha_H>0$ and $\beta_H>0$ are two numerical constants depending only on $H$.

\end{thm}

\medskip

Clearly, the random weight function $\mathcal{Y}$ vanishes, if $(a'\sigma -a \sigma')(x)=0$ for all $x \in \R$. However, if equation (\ref{eq1}) satisfies  
$$ \textrm{(ND)} \qquad \int_{0}^{1}  | \EX \mathcal{Y}_t | \, dt  > 0, $$
then the exact rate of convergence of  the conditional expectation
is $n^{-H-1/2}$, which is summarized in the following Corollary. Note that condition (ND) is satisfied, if e.g. $(a'\sigma -a \sigma')(x_0) \neq 0$.

\medskip

\begin{cor}\label{ndq} If (ND) holds, then there exist constants $C_1=C_1(a,\sigma,x_0,H)>0$ and $C_2=C_2(a,\sigma,x_0,H)>0$ such that
$$      C_1 \cdot n^{-H-1/2} \leq    ( \EX |X_1 -\X_1|^{2} )^{1/2} \leq C_2 \cdot n^{-H-1/2}       $$
for all $n \in \N$.
\end{cor}

\medskip

Consequently, the maximum rate of convergence, which can be obtained by an equidistant discretization of the driving fractional Brownian motion, is $n^{-H-1/2}$  in this case.
Moreover, we can now  characterize the difficulty of equation (\ref{eq1}) in terms of its coefficients:
If the drift- and diffusion coefficient commute, then $X_1$ can be approximated perfectly, see Corollary \ref{perf_app}. Otherwise, if $a$ and $\sigma$ do not commute,  then there are initial values $x_0 \in \R$
such that the exact convergence rate of the conditional expectation is $n^{-H-1/2}$.

Theorem \ref{main_2} fits in the known results for the case $H=1/2$.
In the case of a one-dimensional Stratonovich SDE driven by a Brownian motion $W=(W_t)_{ t \in [0,1]}$, i.e.
$$ dV_t = a(V_t) \, dt + \sigma(V_t) \, d W_t, \qquad V_{0}=v_0\in \R, \quad t \in [0,1], $$
it is well known that
\begin{align}  \lim_{n \rightarrow \infty} \, n \cdot (\EX |V_1- \widehat{V}^{(n)}_1|^{2} )^{1/2} = \frac{1}{\sqrt{12}}  \left(  \int_{0}^{1} \EX | \mathcal{Y}_t^{W} |^{2} \,dt \right)^{1/2},  \label{h_12} \end{align}
where in this case
$$     \widehat{V}^{(n)}_1 = \EX (\, V_1 \, |\,  W_{i/n}, i=0, \ldots, n \, )      $$
and 
$$       \mathcal{Y}_t^{W}= (a\sigma' -a' \sigma) (V_t) \exp \left( \int_{t}^{1} a'(V_\tau) \, d \tau +  \int_{t}^{1} \sigma'(V_{\tau}) \, d W_{\tau}      \right), \qquad t \in [0,1]  .     $$
See, e.g. \cite{cam_clark}, \cite{new1} and \cite{cam_hu}. In these articles, the key for the proof of (\ref{h_12}) is to derive the asymptotic error distribution of $n(V_1- \widehat{V}_1^{(n)})$, which is a conditional normal
distribution with zero mean  and variance  given in terms of  a stochastic differential equation.  Since this method essentially relies on the properties of  Brownian motion,
we could not imitate it and have to rely on indirect methods.

\subsection{Remarks}\label{remarks_section4} 
(i) If the diffusion coefficient is  constant, the boundedness of the drift coefficient is not required in Theorem \ref{main_1} and \ref{main_2}, since we do not have to apply the Lamperti transformation in this case. See also \cite{diss}.
In particular,  we obtain  for  the Langevin equation
$$ dX_t = \lambda X_t \, dt + dB_t, \qquad t \in [0,1], \qquad X_0=x_0$$
that
\begin{align*} \liminf_{n \rightarrow \infty} \,  n^{H+1/2} \,  (\EX |X_1- \X_1|^{2} )^{1/2} 
 \, \geq  \,  \alpha_H \left(\int_{0}^{1} |\mathcal{Y}_t|^{2} \, dt \right)^{1/2} \end{align*}
and
\begin{align*} \limsup_{n \rightarrow \infty} \,  n^{H+1/2} \,  (\EX |X_1- \X_1|^{2} )^{1/2} 
 \, \leq  \,  \beta_H \left(\int_{0}^{1} | \mathcal{Y}_t|^{2} \, dt \right)^{1/2}  \end{align*}
with
$$ \mathcal{Y}_t = \lambda  \exp ( \lambda (1-t)), \qquad t \in [0,1].$$

(ii) Condition (ND) is a non-degeneracy condition for the first term of the chaos expansion of $X_1$. Since $(D_sX_1)_{s \in [0,1]}$ is bounded, we have by (\ref{chaos}) that
$$     X_1 = \EX X_1  + \int_{0}^{1} \EX D_s X_1 \, dB_s + \ldots,  $$
If   (ND) is not satisfied, then we have  
$$          \int_{0}^{1} \EX D_s X_1 \, dB_s = c \cdot B_1 $$
for a constant $c \in \R$, which means that the first chaos of $X_1$ is just a multiple of $B_1$. 
Compare with Theorem \ref{main_1},  which yields, if  $a$ and $\sigma$ commute, the chaos expansion
$$       X_1=    c_0  +  c_1 \cdot B_1 +  \sum_{i=2}^{\infty} c_i \cdot H_i(B_1), $$ 
where $c_i \in \R$ and $H_i$ is the $i$-th Hermite polynomial.

We strongly suppose that the conditions $a'\sigma - a \sigma' \neq 0$ and (ND) are equivalent, i.e. either the error of the conditional expectation is zero or otherwise its exact rate of convergence is $n^{-H-1/2}$.

(iii) At first view, it may seem restrictive to consider only equidistant discretizations of the driving fractional Brownian motion. However, for $H \neq 1/2$ the increments of fractional Brownian motion are  correlated and  therefore  the exact simulation of $B_{t_{1}}, \ldots, B_{t_{n}}$ is in general computationally very expensive. 
Given $n$ iid standard normal random numbers, the Cholesky decomposition method, which is to our best knowledge the only known exact method for the non-equidistant simulation of fractional Brownian motion,   requires still  $O(n^{2})$ arithmetic operations after precomputation of the factorization of the covariance matrix.
However, if  the discretization is equidistant, i.e., $t_{i}=i/n$, $i=1, \ldots, n$,
 the computational cost can be lowered considerably, making use of the stationarity of the increments of fractional Brownian motion. For example, the Davies-Harte algorithm for the equidistant simulation of fractional Brownian motion has 
computational cost $O(n \log(n))$, see e.g. \cite{Craig}. 

For a comprehensive survey of simulation methods for fractional Brownian motion we refer to \cite{Coeur}.

(iv)  For multi-dimensional Stratonovich SDE, i.e.,
$$ dV_t = a(V_t)\, dt  + \sum_{i=1}^{m} \sigma^{(i)}(V_t) \,d W^{(i)}_t $$
with $a,\sigma^{(i)}:  \RR^{d} \rightarrow \RR^{d}$, $i=1, \ldots, m$ and $m$ independent Brownian motions $W^{(i)}$, $i=1, \ldots, m$, 
 the optimal rate of convergence, which can be obtained by point evaluation of the driving Brownian motions, depends on whether the diffusion coefficients $\sigma^{(i)}$ commute or not.
If they commute, that is, if we have
 \begin{align} \label{comm} \sum_{k=1}^{d} \sigma^{(i)}_{k} \frac{\, \, \partial{\sigma^{(j)}_l}} {\partial x_k} = \sum_{k=1}^{d} \sigma^{(j)}_{k} \frac{\, \, \partial{\sigma^{(i)}_l}} {\partial x_k}, \qquad  i,j =1, \ldots, m, \, l=1, \ldots, d, \end{align}
then the optimal rate of convergence is  $n^{-1}$ as in the one-dimensional case, where $n$ is the number of evaluations of the driving multi-dimensional Brownian motion.
(Here $\sigma^{(i)}_k$ denotes the $k$-th component of $\sigma^{(i)}$.)
 However, if (\ref{comm}) is not satisfied, then the optimal rate of convergence is $n^{-1/2}$. See e.g. \cite{cam_clark}.

The prototype example for the latter case is the two-dimensional Stratonovich SDE
\begin{align*}
 dV^{(1)}_t &= V^{(2)}_t \, d W^{(1)}_t , \qquad V^{(1)}_0 =0, \\
 dV^{(2)}_t  &= \, d W^{(2)}_t,   \qquad \qquad V^{(2)}_0 =0.
\end{align*}
Clearly, we have $V_1^{(1)} = \int_{0}^{1} W^{(2)}_t \, dW^{(1)}_t$ and it is well known that
$$  (  \EX | V^{(1)}_1 - \EX (  \,    V^{(1)}_1 \, | \,  W^{(1)}_{i/n}, W^{(2)}_{i/n}, \, i=1, \ldots ,n \,) |^{2} )^{1/2}={\sqrt{1/2}} \cdot n^{-1/2},     $$
 see e.g. \cite{cam_clark}

In a forthcoming paper we will study, whether this  phenomenon also appears for SDEs driven  by fractional Brownian motion.

\subsection{Proof of the lower error bound in Theorem  \ref{main_2}}

The Wiener chaos decomposition described in Subsection \ref{pre_1} will be the key for the proof of the lower error bound, since 
it provides a linearization of the problem, which we have to analyse.

\medskip

\begin{prop}\label{linearize_1}
We have
\begin{align}
\label{linearize}
 \EX | X_1- \X_1|^{2} \geq \EX \left| \int_{0}^{1} \left( B_t - \EX( B_t \, | \, B_{i/n}, i=0, \ldots, n) \right) \EX \mathcal{Y}_t  \,dt   \right|^{2}. \end{align}
\end{prop}

      \medskip

{\noindent \bf Proof.} Clearly, we have that 
$$ \mathcal{J}_{1} (X_1) = \delta( \EX D X_1 ), $$
see equation (\ref{chaos}) in  Subsection \ref{pre_1}.
Recall that 
$$D_s X_1 = \sigma(X_1) \exp \left( \int_{s}^{1} \left(a' -\frac{a \sigma'}{\sigma}\right)(X_{\tau}) \, d \tau \right) $$
by Proposition \ref{rep_mall} and define the mapping $m :[0,1] \rightarrow \RR$ by
$$ m(s)= \EX  D_s X_1, \qquad s \in [0,1]. $$
Since $\sigma$ and $a' -\frac{a \sigma'}{\sigma}$ are bounded due to our assumptions, it follows by dominated convergence 
that the mapping  $m: [0,1] \rightarrow \RR$ is continuously differentiable with
$$ m'(s)= \EX \, D_s X_1 \left(a' - \frac{a \sigma'}{\sigma}\right)(X_s), \qquad s \in [0,1].$$
Using (\ref{mall_deriv}) we obtain
$$ m'(s) = \EX \mathcal{Y}_s, \qquad s \in [0,1].$$
Since the Skorohod integral and the Riemann-Stieltjes integral coincide for smooth  deterministic integrands, see (\ref{rel}), we finally have that
$$   \mathcal{J}_1(X_1 ) =  \int_{0}^{1} \EX \, D_t X_1 \, dB_t = B_1 \EX \sigma(X_1) - \int_{0}^{1} B_t \EX \mathcal{Y}_t  \,dt. $$

The projection of $\X_1 =  \EX( \, X_1 \, | \, B_{i/n}, i=0, \ldots, n \, )$ to the first chaos is  given by the second term of 
the Hermite series expansion of the conditional expectation:
$$ \X_1 = \EX \X_1 + \sum_{i=1}^{n} a_i B_{i/n} + \ldots $$ 
with $a_i \in \RR$ for $i=1, \ldots, n$. See e.g. \cite{p_mall}. Hence we have  
$$  \mathcal{J}_1(\X_1 ) = \sum_{i=1}^{n} a_i B_{i/n}  $$
and  it follows 
$$ \EX | X_1- \X_1|^{2} \geq \EX \left|  B_1 \EX \sigma(X_1) - \int_{0}^{1}  B_t \EX \mathcal{Y}_t  \,dt   - \sum_{i=1}^{n} a_i B_{i/n} \right|^{2}$$
by linearity of the projection to the first chaos. The term on the right hand side is the
error of the quadrature formula $$ \widehat{I}_n(B)=B_1 \EX \sigma(X_1) - \sum_{i=1}^{n} a_i B_{i/n}$$  for the approximation of the integral $$ I(B)=\int_{0}^{1}  B_t \EX \mathcal{Y}_t  \,dt.$$
Clearly, by definition of the conditional expectation the best quadrature formula, which uses $B_0, B_{1/n}, \ldots, B_1$,  is given
by $$\overline{I}_n(B) = \int_{0}^{1}  \EX( B_t \, | \, B_{i/n}, i=0, \ldots, n) \,   \EX \mathcal{Y}_t \, dt.$$
Thus, we finally obtain
$$ \EX | X_1- \X_1|^{2} \geq \EX \left| \int_{0}^{1}  B_t \EX \mathcal{Y}_t  \,dt -  \int_{0}^{1}  \EX( B_t \, | \, B_{i/n}, i=0, \ldots, n) \,   \EX \mathcal{Y}_t \, dt   \right|^{2}.$$
\fin

To finish the proof of the lower error bound in Theorem \ref{main_2}, we  now have to analyse the quantity
$$ e(n)= \left( \EX \left| \int_{0}^{1} \EX \mathcal{Y}_s (B_s  - \EX( B_s | B_{i/n}, i=0, \ldots, n)) ds \right|^{2} \right)^{1/2}.$$
For this,  we recall some well known facts about reproducing Kernel Hilbert spaces and weighted integration problems for Gaussian processes. See e.g.~ \cite{ritter} and the references therein.

Let $\mathcal{W}=(\mathcal{W}_t)_{t \in \R}$ be a continuous stochastic process with covariance function $(K(s,t))_{s,t \in \R}$.  Then the reproducing Kernel Hilbert space $(\mathfrak{H}(K), \langle \cdot, \cdot \rangle_{K})$ corresponding to the process $\mathcal{W}$  is the uniquely determined Hilbert space of real valued functions on $\R$ such that
$$ K(\cdot,t) \in \mathfrak{H}(K)$$ 
and
$$ \langle h, K(\cdot,t) \rangle_{K} =h(t), \qquad t \in \R$$
holds for all $h \in \mathfrak{H}(K)$.

Now consider the linear functional $$I(\mathcal{W})= \int_{\R} \rho(t) \mathcal{W}_t \, dt,$$
where the function $\rho \in C(\R)$ is non-negative and has compact support. Then the error of a quadrature formula 
$$\widehat{I}(\mathcal{W})= \sum_{i=1}^{n} a_i \mathcal{W}_{t_{i}} $$
with $a_i, t_{i} \in \R $ for $i=1, \ldots, n$, can be characterized  in terms of the  reproducing Kernel Hilbert space as follows:
$$ (\EX | I(\mathcal{W})-\widehat{I}(\mathcal{W} ) |^{2} )^{1/2} = \sup \left \{ \left| \int_{\R} \rho(t)h(t) dt   -\sum_{i=1}^{n}a_{i} h(t_{i}) \right| ; \,   h \in \mathfrak{H}(K): \, \| h\|_{K} \leq 1 \right \}. $$
In addition, if the process $\mathcal{W}$ is Gaussian and if we use the conditional expectation of the functional given the evaluations of the random process, i.e.
$$ \overline{I}(\mathcal{W})=  \EX (\,  I(\mathcal{W})   \, | \, \mathcal{W}_{t_{i}}, i=1, \ldots, n \,), $$
as quadrature formula, then we have
\begin{align*}  &(\EX |I(\mathcal{W})-\overline{I}(\mathcal{W})^{2} )^{1/2} \\ &\qquad = \sup \left \{  \left| \int_{\R} \rho(t)h(t) \, dt \right|; \, h \in \mathfrak{H}(K): \, \| h\|_{K} \leq 1, \, h(t_{i})=0, \, i=1, \ldots, n \right \} .\end{align*}
See e.g.~section III.2 in \cite{ritter}. Now we can finish the proof of (\ref{low_bound_formula}):

{ \noindent \bf Proof of the lower bound in Theorem \ref{main_2}:} By the above considerations it remains to analyse the quantity
$$ e(n)=  \sup \left \{  \left| \int_{0}^{1} \rho(t)h(t) \, dt\right| ; \, h \in \mathfrak{H}(R_H): \, \| h\|_{R_H} \leq 1, \, h(i/n)=0, \, i=0, \ldots, n \right \}, $$
with $\rho(t) = \EX \mathcal{Y}_t$, $t \in [0,1]$. Here 
$ (\mathfrak{H}(R_H), \langle \cdot, \cdot \rangle_{R_H})$ is the reproducing kernel Hilbert space of  $B=(B_t)_{t \in \R}$, which is the fractional Brownian motion defined on the whole real line.

It is well known that
$$ G=\{g \in C^{\infty}(\R):  \, \textrm{supp} \,g \, \textrm{ compact}, \, \, 0 \notin \textrm{supp} \,g  \} $$
is dense in $\mathfrak{H}(R_H)$ and that
$$ \|h \|_{R_H}^{2}= c_H \int_{\R} | u |^{2H+1} | \mathcal{F}(h)(u) |^{2} \, du $$
for $h \in G$, where the constant $c_H >0$ is known explicitly, see e.g.~ section 6.1 in \cite{ritter} and \cite{singer}. Here  $\mathcal{F}(h)$ denotes the Fourier transform of $h$, defined by
 $$ \mathcal{F}(h)(t) = \int_{\R} \exp(- i tu) h(u) \, du. $$

Now consider the stationary Gaussian process $ \mathcal{V}=(\mathcal{V}_t)_{t \in \R}$ with  covariance kernel
$$K(0,t)= \int_{\R} \exp( i tu) f(u) \, du, \quad t \in\R, $$
where the spectral density $f: \R \rightarrow \R$ is given by
\begin{eqnarray} \label{spec_dens} f(u)= \frac{1}{4\pi^2 c_H}(1+u^{2})^{-H-1/2}, \qquad u \in \R. \end{eqnarray}
Then the reproducing Kernel Hilbert space corresponding to $\mathcal{V}$ is the Bessel potential space $$\mathfrak{B}^{H+1/2}=\{ h \in L^{2}(\R): |h|_{H+1/2} < \infty \}$$
with
$$ |h|_{H}^{2}= \frac{1}{4 \pi^2} ¸\int_{\R} f(u)^{-1} |\mathcal{F}(u)|^{2} \, du. $$

We have $$ \{ h \in \mathfrak{B}_{H}: h(0)=0\} \subset \mathfrak{H}(R_H),$$
and 
 $$ \| h \|_{R_H} \leq |h|_{H} $$ 
for $h \in \mathfrak{B}_{H}$, see e.g.~ section 6.1 in \cite{ritter}.
Hence it follows $$ e(n) \geq  \sup \left \{  \left| \int_{0}^{1} \rho(t)h(t) \, dt \right| ; \,  h \in \mathfrak{B}_{H}: \, \| h\|_{R_H} \leq 1, \, h(i/n)=0, \, i=0, \ldots, n \right \}$$
and moreover
$$ e(n) \geq  \sup \left \{  \left| \int_{0}^{1} \rho(t)h(t) \, dt \right| ; \,  h \in \mathfrak{B}_{H}: \, | h|_{H} \leq 1, \, h(i/n)=0, \, i=0, \ldots, n \right \}.$$

However, the quantity on the right hand side of the above  equation corresponds to
$$ \widetilde{e}(n)=\left( \EX \left| \int_{0}^{1} \rho(t) (\mathcal{V}_t - \EX(\mathcal{V}_t \, | \, \mathcal{V}_{i/n}, i=0, \ldots, n)) \,dt \right|^{2} \right)^{1/2},$$
where $\mathcal{V}$ is the stationary Gaussian process with covariance kernel given by (\ref{spec_dens}).

From \cite{stein}, Proposition 2.1,  we now obtain that
$$ \lim_{n \rightarrow \infty} \frac{\widetilde{e}(n)}{s(n)}=1 $$
with
$$s(n)=\int_{-n\pi}^{n\pi} \sum_{j \in \mathbb{Z} \setminus \{0 \} } f(u+n 2\pi j) |\mathcal{F}(\rho(u))|^{2} \, du.$$
We have
$$s(n) \geq   \int_{-n\pi}^{n\pi} |\mathcal{F}(\rho(u))|^{2} \, du \cdot \sum_{j \in \mathbb{Z} \setminus \{ 0 \} } f(\xi+n 2\pi j)$$
with $\xi \in (-n\pi, n \pi)$. Note that by symmetry of the spectral density it holds
$$ \sum_{j \in \mathbb{Z} \setminus \{0 \} } f(u+n 2\pi j) | = 2 \sum_{j=1}^{\infty} f(u+n 2\pi j), \qquad u \in \RR. $$
Since $f$ is strictly decreasing, it follows
$$s(n)\geq  2\sum_{j=1}^{\infty} f(n (2\pi j +\pi)) \int_{-n\pi}^{n\pi}  |\mathcal{F}(\rho(u))|^{2} \, du.$$

Now consider the quantity $n^{2H+1}s(n)$. Clearly, we have
$$  \liminf_{n \rightarrow \infty}n^{2H+1}s(n) \geq  \frac{1}{2\pi^{2}c_H}  \liminf_{n \rightarrow \infty} \sum_{j=1}^{\infty} (n^{-2}+ \pi(2j+1)^{2})^{-H-1/2} \int_{\R} |\mathcal{F}(\rho(u))|^{2} \, du. $$
Fatou's Lemma and Parseval's equality yield
$$  \liminf_{n \rightarrow \infty}n^{2H+1}s(n) \geq  \frac{1}{c_H 2 \pi^{2}}  \sum_{j=1}^{\infty} \frac{\pi}{((2j+1)\pi )^{2H+1}}  \int_{\R} |\rho(u)|^{2} \, du. $$

Since
$$  \liminf_{n \rightarrow \infty}n^{2H+1}e(n) \geq \liminf_{n \rightarrow \infty}n^{2H+1}s(n),$$
it  follows 
$$ \liminf_{n \rightarrow \infty}n^{2H+1}e(n) \geq 
 \alpha_H^{2} \int_{\R} |\rho(u)|^{2} \, du$$
with
$$ \alpha_H^{2} = \frac{1}{2\pi^{2H+2}} \frac{1}{c_H}  \sum_{j=1}^{\infty}\frac{1}{(2j+1)^{2H+1}}.$$
Since $\rho(t)= \EX \mathcal{Y}_t$, $t \in [0,1]$, and
$$ e(n) \leq \EX |X_1- \X_1|^{2}$$ by Proposition \ref{linearize_1}, we finally have shown  the lower bound in Theorem \ref{main_2}.

\fin

\subsection{Proof of the upper bound in Theorem \ref{main_2} }
For the proof of (\ref{up_bound_formula}) we will again use an indirect method.

Denote by $\BL=(\BL_t)_{t \in [0,1]}$  the piecewise
linear interpolation of $B$ based on an equidistant discretization with stepsize $1/n$, that is
$$ \BL_t = B_{k/n} +  (  nt  -k )   \left( B_{(k+1)/n} - B_{k/n} \right), \qquad t \in [k/n , (k+1)/n],$$
and consider the following stochastic differential equation
\begin{align*} 
 d \XL_t&= a(\XL_t) \, dt + \sigma(\XL_t) \, d\BL_t, \qquad t \in [0,1], \\ \XL_0&=x_{0}. \nonumber
\end{align*}
This equation has a pathwise unique solution, since $\BL$ has finite total variation for every $n \in \N$.
In particular, we can again apply the Lamperti transformation given in Subsection \ref{pre_2}
and obtain that
$$ \XL_t = \vartheta^{-1}(\YL_t), \qquad t \in [0,1],$$
where $\YL=(\YL_t)_{t \in [0,1]}$ is the unique solution of
\begin{align*}
d\YL_t=& g(\YL_t) \, dt + d\BL_t, \qquad t \in [0,1], \\
\YL_0 =&  \vartheta(x_{0}).
\end{align*}

Clearly, we have
$$ \EX |X_1 -\X_1|^{2} \leq \EX |X_1 -\XL_1|^{2} $$
for every $n \in \N$. The next Proposition characterizes the leading term of the error of  the approximation $\XL_1$:

\medskip

\begin{prop}\label{linearize_2}
We have
\begin{align} \label{linearize_2eq}
\EX |X_1 -\XL_1|^{2} = \EX \left| \int_{0}^{1}  \mathcal{Y}_t (B_t- \BL_t) \, dt \right|^{2} + O(n^{-2H}).
\end{align}
\end{prop}

\medskip

{ \noindent \bf Proof.}  We will denote constants, which depend only on $x_0$, $H$, $a$, $\sigma$ and their derivatives by $c$, regardless of their value.

(i) We  first we establish the following estimate:
\begin{eqnarray} \label{bound_Y}
\sup_{t \in [0,1]} |Y_t|   + \sup_{t \in [0,1]} |\YL_t| \leq c \sup_{t \in [0,1]} |B_t|.
\end{eqnarray}
For this note that
$$ |Y_t| \leq  |\vartheta^{-1}(x_0)| + \int_{0}^{t} g(Y_\tau) \, d \tau + |B_t|, \qquad t \in [0,1]. $$
Since $g$ is continuously differentiable with bounded derivative, it satisfies a linear growth condition
and we obtain
$$ \sup_{s \in [0,t]} |Y_s| \leq c  + c \int_{0}^{t} \sup_{s \in [0,\tau]} |Y_s| \,  d \tau +  \sup_{s \in [0,t]} |B_s|, \qquad t \in [0,1].$$ 
Thus, Gronwall's Lemma yields that  
$$ \sup_{t \in [0,1]} |Y_t| \leq c \sup_{t \in [0,1]} |B_t|.$$
Similar, we obtain
 $$ \sup_{t \in [0,1]} |\YL_t| \leq c \sup_{t \in [0,1]} |B_t|,$$
since
$$\sup_{t \in [0,1]} |\BL_t| \leq \sup_{t \in [0,1]} |B_t|.$$

(ii) Using the Lamperti transformation, we can write
$$ X_t- \XL_t = \vartheta^{-1}(Y_t) - \vartheta^{-1}(\YL_t), \quad t \in [0,1].$$
Since we have
\begin{align*}
(\vartheta^{-1})'(x)&=\sigma(\vartheta^{-1}(x)),\\
(\vartheta^{-1})''(x)&= \sigma \sigma'(\vartheta^{-1}(x)) , 
\end{align*}
 we  obtain
\begin{eqnarray*} 
X_t- \XL_t =\sigma(X_t)(Y_t-\YL_t) + \frac{1}{2}\sigma\sigma'(\vartheta^{-1}( \theta_t+ (1-\theta_t) \YL_t))(Y_t-\YL_t)^{2}
\nonumber \end{eqnarray*}
for $t\in [0,1]$ with a random $\theta_t \in (0,1)$.
So we only need to consider the difference between $Y$ and $\YL$.
For $$Z_t= Y_t - \YL_t, \qquad t \in [0,1],$$ we obtain
\begin{eqnarray*}
Z_t=    B_t-\BL_t + \int_{0}^{t} g(Y_{\tau})- g(\YL_{\tau}) \, d \tau, \qquad t \in [0,1].
\end{eqnarray*}
The Lipschitz continuity of  $g$ implies that
$$   |Z_t| \leq    |B_t-\BL_t| +  c \int_{0}^{t} |Z_{\tau}| \, d \tau , \qquad t \in [0,1]  $$
and Gronwall's Lemma yields
$$|Z_t| \leq  c \int_{0}^{t} |B_{\tau}- \BL_{\tau}| \, d \tau,   \qquad t \in [0,1].  $$
Thus, it follows 
\begin{align} \label{est_Z}  \sup_{t \in [0,1]} \EX \left| Z_t\right|^p \leq c \cdot n^{-Hp},  \end{align}
since clearly
$$          \sup_{t \in [0,1]} \EX \,  | B_t- \BL_t |^p \leq c \cdot n^{-Hp}.$$
Using the boundedness of $\sigma$ and $\sigma'$ we obtain 
\begin{align} \label{expansion_X} X_t- \XL_t =\sigma(X_t)Z_t  + R^{(1)}_t  \end{align}
with
$$ \sup_{t \in [0,1]} \EX | R^{(1)}_t|^p \leq c \cdot n^{-2Hp}. $$ 

(iii) Now, we analyse the process $Z=(Z_t)_{t \in [0,1]}$ in more detail. We can write
\begin{eqnarray} \label{rec_Z}
Z_t=    B_t-\BL_t + \int_{0}^{t} g'(Y_{\tau})Z_{\tau} \, d \tau +R^{(2)}_{t}, \qquad t \in [0,1],
\end{eqnarray} 
where
$$R^{(2)}_t= \int_{0}^{t}    ( g'(\theta_{\tau} Y_{\tau}+ (1-\theta_{\tau}) \YL_{\tau})    -g'(Y_{\tau}) ) Z_{\tau} \, d \tau, \qquad t \in [0,1], $$ 
with $\theta_{\tau} \in (0,1)$.
Since $g'$ is Lipschitz continuous, it follows
$$     \left|   g'(\theta_{\tau} Y_{\tau}+ (1-\theta)_{\tau} \YL_{\tau})    -g'(Y_{\tau}) \right| \leq c |Z_{\tau}| , \qquad \tau \in [0,1] , $$
and consequently
$$    |R^{(2)}_t| \leq      c \int_{0}^{t}    |Z_{\tau}|^{2} \, d \tau,  \qquad t \in [0,1],  $$
which in turn yields together with (\ref{est_Z}) that
\begin{align}\label{est_r2}
  \sup_{t \in [0,1]}  \EX   |R^{(2)}_t|^p \leq      c  \cdot n^{-2Hp}. \end{align}

Now consider the equation
$$ \widetilde{Z}_t = B_t-\BL_t + \int_{0}^{t} g'(Y_{\tau})\widetilde{Z}_{\tau} \, d \tau , \qquad t \in [0,1].$$
Its unique solution  is given by
\begin{align}\label{sol_diff}
\widetilde{Z}_t = \int_0^t g'(Y_{\tau}) \exp\left( \int_{\tau}^{t} g'(Y_s) \,ds \right)(B_{\tau}-\BL_{\tau}) \, d \tau, \qquad t\in [0,1]. \end{align}
Since \begin{eqnarray*}
Z_t-\widetilde{Z}_t = R^{(2)}_{t} + \int_{0}^{t} g'(Y_{\tau})(Z_{\tau}-\widetilde{Z}_{\tau}) \, d \tau, 
\end{eqnarray*}
again an application of Gronwall's Lemma and (\ref{est_r2}) yield that
\begin{align}\label{est_r3}  \sup_{t \in [0,1]} \EX |Z_t - \widetilde{Z}_t|^{p} \leq c \cdot n^{-2Hp}.\end{align}
Therefore, we obtain from (\ref{sol_diff}) and (\ref{est_r3}) that 
$$ Z_t =  \int_0^t g'(Y_{\tau}) \exp\left( \int_{\tau}^{t} g'(Y_s) \,ds \right)(B_{\tau}-\BL_{\tau}) \, d \tau + R_{t}^{(3)}$$
with
$$    \sup_{t \in [0,1]}        \EX | R^{(3)}_t  |^p \leq c \cdot n^{-2Hp}. $$

(iv) Using (\ref{expansion_X}) we finally have  that
$$       X_t - \XL_t = \sigma(X_t)Z_t +  R_t^{(4)}, \qquad t \in [0,1],   $$
with
$$         \sup_{t \in [0,1]}                            \EX | R^{(4)}_t  |^p \leq c \cdot n^{-2Hp}. $$
Since 
$$ g'(Y_s) = \left(a'- \frac{a\sigma'}{\sigma} \right)(X_s), \qquad s \in [0,1],$$
we obtain
\begin{align*}      X_t - \XL_t  &= \int_0^t     \left( a' - \frac{a\sigma'}{\sigma} \right)(X_\tau)     \sigma(X_t) \\ & \qquad \qquad \quad \quad \times \exp \left( \int_{\tau}^{t}   \left( a' - \frac{a\sigma'}{\sigma} \right)(X_s)    \,ds \right)(B_{\tau}-\BL_{\tau}) \, d \tau +  R_t^{(4)}    \end{align*}
for $t \in [0,1]$,
and the assertion follows, since
\begin{align*}
   &  \sigma(X_t) \exp \left( \int_{\tau}^{t}   \left( a' - \frac{a\sigma'}{\sigma} \right)(X_s)    \,ds \right) 1_{[0,t]}(s) \\ & \qquad \qquad \qquad \qquad =       \sigma(X_{\tau}) \exp \left( \int_{\tau}^{t}   a'(X_s)    \,ds  + \int_{\tau}^{t}  \sigma'(X_s) \, d B_s \right)         1_{[0,t]}(s)    \end{align*}
by  (\ref{rep_mall_formula}).
\fin

Now, it remains to analyse the quantity
$$    \EX \left| \int_{0}^{1}  \mathcal{Y}_t (B_t- \BL_t) \, dt \right|^{2}.    $$
Note that this is the mean square error of the trapezoidal-type quadrature formula
\begin{align*}\widetilde{I}(B)& =  n \int_{0}^{t_1}  (t_1-t) \mathcal{Y}_t \, dt \cdot    B_0 \\ & \qquad  \qquad +    n \sum_{i=1}^{n-1}      \left(   \int_{t_i}^{t_{i+1}}  (t_{i+1}-t) \mathcal{Y}_t \, dt  +    \int_{t_{i-1}}^{t_{i}}  (t-t_{i-1}) \mathcal{Y}_t \, dt     \right) \cdot B_{i/n}     \\ & \qquad  \qquad \qquad \quad    +   n \int_{t_{n-1}}^{1} (t-t_{n-1}) \mathcal{Y}_t \, dt \cdot B_1,   \end{align*}
where $t_i =i/n$, $i=0, \ldots ,n$,
for the approximation of the weighted integral 
$$ I(B)= \int_{0}^{1} \mathcal{Y}_t B_t \,dt $$ with
random weight function $\mathcal{Y}=(\mathcal{Y}_{t})_{ t \in [0,1]}$.

Combining the  next Proposition and Proposition \ref{linearize_2}   we obtain the upper bound of Theorem \ref{main_2}.

\medskip

\begin{prop} \label{int_problem} We have
  $$  \lim_{n \rightarrow \infty} \, n^{2H+1} \EX \left| \int_{0}^{1}  \mathcal{Y}_t (B_t- \BL_t) \, ds \right|^{2} = |\zeta(-2H)|    \int_{0}^{1} \EX | \mathcal{Y}_t |^{2}  dt,  $$
where $\zeta$ denotes the Riemann Zeta function.
\end{prop}

\medskip

The proof of this Proposition is rather technical. Therefore we postpone it to the Appendix.
Note that for the constant $\beta_H$ in Theorem \ref{main_2} we have
$$ \beta_H^{2} = |\zeta(-2H)|. $$

\section{McShane's scheme method}\label{sec_mac}
Theorem \ref{main_2} states  in particular that  the  rate of convergence of the
conditional expectation $\X_1$ is  at least $n^{-H-1/2}$. 
 Clearly, the conditional expectation is explicitly known only in some exceptional
cases and   thus  is not an implementable approximation scheme in general .
In this section, we present a  feasible approximation scheme, which is almost
as good as the conditional expectation in the sense that its convergence rate is also at least $n^{-H-1/2}$.

The McShane's scheme for stochastic differential equation driven by fractional Brownian motion is defined by
$\XM_{0}=x_{0}$ and \begin{align} \label{ms} \XM_{k+1}=\XM_{k} &+ a(\XM_{k}) \Delta + \sigma(\XM_{k})\Delta_k B + \frac{1}{2} \sigma \sigma'(\XM_k) (\Delta_k B)^{2} \\ &+\frac{1}{2}(a\sigma' + a'\sigma)(\XM_k) \Delta_k B \Delta + \frac{1}{2} aa'(\XM_k) \Delta^{2} \nonumber  \\ & + \frac{1}{6} (\sigma^{2} \sigma'' +  \sigma(\sigma')^{2})(\XM_k)(\Delta_k B)^{3}  \nonumber \end{align} for $k=0, \ldots, n-1$, where $\Delta=1/n$ and  $\Delta_k B = B_{(k+1)/n} - B_{k /n}$.   
For Stratonovich SDEs driven by Brownian motion, this scheme was studied e.g.~
in \cite{macshane}, \cite{new1} and shown to be asymptotically efficient,  i.e.
\begin{align*} 
\lim_{ n \rightarrow \infty} n \,  \EX | U_1 - \overline{U}_1^{(n)}  |^{2} =
\lim_{ n \rightarrow \infty} n \,  \EX | U_1 -\widehat{U}^{(n)}_n |^{2} = \frac{1}{12}
\int_{0}^{1} \EX | \mathcal{Y}_t^{W}|^{2} \, dt .
\end{align*}

The following Theorem gives  in the non-degenerated case the exact convergence rate for McShane's method  for
SDEs driven by fractional Brownian motion:

\medskip

\begin{thm} \label{main_thm_3} For the approximation scheme given by (\ref{ms}) 
  we have
\begin{align}
\lim_{n \rightarrow \infty} n^{H+1/2} (\EX | X_1 - \XM_n|^2)^{1/2}  =  \beta_H
 \left(  \int_{0}^{1} \EX | \mathcal{Y}_t|^{2} \, dt
\right)^{1/2}. \label{opt_rate_cons}
\end{align} 
\end{thm}

\medskip

Note that the asymptotic constant on the right hand side of  (\ref{opt_rate_cons}) vanishes, if and only if $(a'\sigma - a \sigma')(x)=0$ for all $x \in \R$. See Remark (iii) in Section \ref{remarks_section5}.
Hence, if condition (ND) holds, then both McShane's method and the conditional expectation have exact rate of convergence $n^{-H-1/2}$ and thus McShane's method is optimal in this case.

In the degenerated case, we obtain the following upper bound for the
error of McShane's scheme.

\medskip

\begin{prop} If $(a'\sigma - a \sigma')(x)=0$ for all $x \in \RR$, then there exists a constant $C=C(a,\sigma,x_0,H)>0$ such that
$$  ( \EX |X_1 - \XM_n |^{2} )^{1/2}  \leq C \cdot n^{-2H}$$
for all $n \in \N$.
\end{prop}

\medskip

{ \noindent \bf Proof.}  This is a straightforward consequence of Proposition \ref{linearize_2}.
\fin

Thus, in this case McShane's scheme has a convergence rate of at least $n^{-2H}$.  For the case of a zero drift coefficient, i.e. $a=0$, it is shown in \cite{MG} that
$$       n^{4H-1} \left[ X_1 - \XM_n \right] \stackrel{Prob}{\longrightarrow}  - \frac{1}{8} \sigma(X_1) \int_{0}^{1} ( \sigma'^{3} + 4\sigma \sigma' \sigma'' + \sigma^{2}\sigma''')(X_s) \, ds   $$ 
for $n \rightarrow \infty$.
 So we strongly suppose that the exact convergence rate of McShane's scheme in the degenerated case is $n^{-4H+1}$.

\subsection{Remarks} \label{remarks_section5} 

(i) In this article, we use indirect methods to  determine
 the mean-square error of  the considered approximation schemes.  The main reason for this is  that moment estimates for Riemann-Stieltjes integrals driven by fractional Brownian motion are much more involved than for It\^{o} integrals with respect to Brownian motion.

Recall that, if the process $u=(u_t)_{t \in [0,1]}$ satisfies appropriate smoothness conditions, the relation between the Riemann-Stieltjes integral and the Skorohod integral is given by
\begin{align}\label{rel_ex} \int_{0}^{1} u_t \, dB_t = \delta(u) + \alpha_H \int_{0}^{1} \int_{0}^{1} D_s u_t |s-t|^{2H-2} \, ds \, dt \end{align}
and moreover 
\begin{align}
\label{iso_work_ex}\EX \left|  \delta (u) \right|^{2} &= \, \alpha_{H} \int_{[0,1]^{2}} \EX \,u_s u_r |s-r|^{2H-2} \, dr \, ds \\ &
 \qquad  + \alpha_{H}^{2} \int_{[0,1]^{4}}  \EX D_{r}u_s D_{r'}u_{s'} |r-s'|^{2H-2} |r'-s|^{2H-2} dr dr' ds ds'
\nonumber.\end{align}

Since in both expressions the Malliavin derivative appears, it is not possible to use them for a direct error analysis. To illustrate this, consider e.g. the continuous  Euler method 
for equation (\ref{eq1}), 
which is given by $X_0^{E}=x_0$ and
$$  X_t^{E} = X_{k/n}^{E} +    \int_{k/n}^{t} a( X_{k/n}^{E} ) \, dt   + \int_{k/n}^{t} \sigma( X_{k/n}^{E} ) \, dB_t ,\qquad t \in [k/n, (k+1)/N)$$
for $k=0, \ldots , n-1$. Here we have
$$     X_t - X_t^{E} =     \int_{0}^{t}    a(X_t)-   a( X_{[nt]/n}^{E} ) \, dt +   \int_{0}^{t}    \sigma(X_t)-   \sigma( X_{[nt]/n}^{E} ) \, dB_t, \qquad t \in [0,1] .  $$
Applying (\ref{rel_ex}) and (\ref{iso_work_ex}) would yield an equation, which 
involves the first Malliavin derivative of $X_t - X_t^{E}$ for $t \in [0,1]$. Thus, for analyzing the Euler method in this way, we would need  to control the difference between the Malliavin derivative of the solution and  the Malliavin derivative of the Euler method. But this  involves the
second Malliavin derivative etc.~ and we cannot have closable formulas.

In this article. we apply  the Lamperti transformation to avoid this problem. Essentially the Lamperti transformation reduces the error analysis of the considered equation to the error analysis of a
related
 equation with constant diffusion coefficient, for which the above problem does not appear. The price we have to pay for this procedure is the quite strong assumption (A3)  on the drift- and diffusion coefficient, which ensures the integrability of the remainder terms in the error analysis.

If the pathwise error is considered instead of the mean square error, then assumption (A3) can be avoided. 
Under the assumptions (A1) and (A2) and using the Doss-Sussmann transformation as reduction method to analyse the Euler scheme, it is shown in \cite{NN}  that
$$      n^{2H-1} (  X_1 - X_1^{E}  ) \stackrel{a.s}{\longrightarrow} - \frac{1}{2} \int_{0}^{1} \sigma'(X_s) D_s X_1 \, ds        $$
for $n \rightarrow \infty$.

(ii) For the proof of  Theorem \ref{main_thm_3} we will use Proposition \ref{linearize_2},
i.e.  
$$ \EX |X_1 -\XL_1|^{2} = \EX \left| \int_{0}^{1}  \mathcal{Y}_t (B_t- \BL_t) \, dt \right|^{2} + O(n^{-2H}),$$
where
\begin{align} \label{lin_eq} 
 d \XL_t&= a(\XL_t) \, dt + \sigma(\XL_t) \, d\BL_t, \qquad t \in [0,1], \\ \XL_0&=x_{0}, \nonumber
\end{align}
and $\BL$ is  the piecewise linear interpolation of $B$ with stepsize $1/n$.
Therefore, it remains to compare $\XL_1$ and $\XM_n$ to show   Theorem \ref{main_thm_3}.
Note that equation (\ref{lin_eq}) actually corresponds to a system of piecewise random ordinary differential equations: If we define
\begin{align*} 
\dot{x}^{(k)}(t) &= a(x^{(n)}(t)) +     n( B_{(k+1)/n} -B_{k/n})  \sigma(x^{(n)}(t)), \qquad t \in (k/n, (k+1)/n]  \\
      x^{(k)}( k/n) &= x^{(k-1)}(k/n)
\end{align*}
for $k=0, \ldots n-1$ with $x^{(0)}(0) = x_0$, then we have
$$         \XL_t = x^{(k)}(t) , \qquad t \in [k/n, (k+1)/n]. $$
Hence Theorem \ref{main_thm_3} will also hold for any method, which approximates the above system of piecewise ordinary differential equations with convergence rate $n^{-H-1/2- \varepsilon}$ for $\varepsilon >0$.  Compare with  \cite{cas_gaines} for the case $H=1/2$.

(iii) Since  $$ \mathcal{Y}_t = (a\sigma' - a'\sigma)(X_t) \exp \left( \int_{t}^{1} a'(X_s) \,ds 
+ \int_{t}^{1} \sigma'(X_s) \, dB_s \right), \qquad t \in [0,1],$$
the condition  $ (a\sigma' - a'\sigma)(x)=0$ for all $x \in \R$ clearly implies $ \int_{0}^{1} \EX |\mathcal{Y}_t|^{2} \, dt =0$.

On the other hand, assume that  $ \int_{0}^{1} \EX |\mathcal{Y}_t|^{2} \, dt =0$. However,  Proposition \ref{rep_mall} now yields that
$$\int_{0}^{1} \EX | (a\sigma' - a'\sigma)(X_t)      |^{2} \, dt =0.$$
Since $\sigma$ is strictly positive,  the distribution of  $X_t$ is absolutely continuous with respect to the Lebesgue measure with a strictly positive density for every $t \in (0,1]$. See e.g \cite{nosi}. Hence  it follows $ (a\sigma' - a'\sigma)(x)=0$ for all $x \in \R$.

Consequently, we have   $ \int_{0}^{1} \EX |\mathcal{Y}_t|^{2} \, dt =0$ if and only if
$ (a\sigma' - a'\sigma)(x)=0$ for all $x \in \R$.

\subsection{Proof of the convergence rate of McShane's method in the non-degenerated case}

To show (\ref{opt_rate_cons}), we require the following Lemmata, which
both can be shown  by straightforward calculations.
The first one  gives a
Taylor expansion of the Lamperti transformation, while the second considers
the regularity of the drift coefficient of the reduced equation:

\medskip

\begin{lem}
For the mapping $\vartheta: \R \rightarrow \R$  given by 
$$ \vartheta(x) = \int_{0}^{x} \frac{1}{\sigma(\xi)} \, d \xi, \qquad x
\in \R,$$
we have $\vartheta \in C^{4}(\R;\R)$ with bounded derivatives. In particular,
it holds
\begin{align*}
\vartheta'(x)&=\frac{1}{\sigma(x)},\\
\vartheta''(x)&= -\frac{\sigma'(x)}{\sigma^{2}(x)}, \\
\vartheta'''(x)&= -\frac{\sigma''(x)}{ \sigma^{2}(x)}+ \frac{2(\sigma')^{2}(x)}{\sigma^{3}(x)}
\end{align*} for $x \in \R$.
\end{lem}

\medskip

\begin{lem} \label{lemma_g} Define $g: \R \rightarrow \R$ by 
$$ g(x)=  \frac{a}{\sigma}(\vartheta^{-1}(x)), \qquad x \in \RR.$$
 Then $g$ is bounded and  we have $g \in C^{2}(\R; \R)$ with bounded derivatives. Moreover, it holds
$$  g'(x)=     \left( a - \frac{a \sigma'}{\sigma}\right)(\vartheta^{-1}(x)), \qquad x \in \R.$$
\end{lem}

\medskip

{\bf \noindent  Proof of Theorem  \ref{main_thm_3}.}  We will denote constants, which depend only on $x_0$, $H$, $a$, $\sigma$ and their derivatives by $c$, regardless of their value.
Using Proposition \ref{linearize_2}, the assertion of  Theorem  \ref{main_thm_3} follows, if we show that
\begin{align}  \EX| \XL_1 - \XM_1 |^{2} \leq c \cdot  \Delta^{4H}
  \label{to_show_cons} .\end{align}

(i) For this, define
$$ \YM_{k}=\vartheta(\XM_k),\qquad k =0, \ldots, n,$$
which turns out to be an approximation scheme for $\YL$: A Taylor
expansion yields
\begin{align*}
\YM_{k+1}= \YM_k &+ \vartheta(\XM_{k+1}) - \vartheta(\XM_k) \\
         = \YM_k &+ \vartheta'(\XM_k)  (\XM_{k+1}-\XM_k)  \\ &+ \frac{1}{2}
         \vartheta''(\XM_k)  (\XM_{k+1}-\XM_k)^{2}  \\ & + 
\frac{1}{6} \vartheta'''(\XM_k)  (\XM_{k+1}-\XM_k) ^{3} \\ &+ \frac{1}{24}\vartheta^{\textrm{(iv)}}(\theta_k \XM_k + (1-\theta_k) \XM_{k+1})  (\XM_{k+1}-\XM_k) ^{4}
\end{align*} for $k=0, \ldots, n-1$
with a random $\theta_k \in (0,1)$ .
Since $a$ and $\sigma$ are bounded together with their derivatives, we moreover
have that
\begin{align} \label{smooth_inc}  | \XM_{k+1}-\XM_k | \leq c \cdot \left(
  \Delta +
  |\Delta_k B| + |\Delta_k B|^{2} + |\Delta_k B|^{3} 
        \right).   \end{align}
Hence we obtain that
\begin{align}
\YM_{k+1} \label{rec_ym}
         = \YM_k &+ \vartheta'(\XM_k)  (\XM_{k+1}-\XM_k)  \\ &+ \frac{1}{2} \nonumber 
         \vartheta''(\XM_k)  (\XM_{k+1}-\XM_k)^{2}  \\ & +  \nonumber
\frac{1}{6} \vartheta'''(\XM_k)  (\XM_{k+1}-\XM_k) ^{3} + R_k^{(1)}
\end{align}
for $k=0, \ldots, n-1$ with 
$$ \sup_{k=0, \ldots, n-1} \EX | R_k^{(1)}|^{p} \leq c \cdot \Delta^{4pH} $$
by the boundedness of $\vartheta^{\textrm{(iv)}}$ and the estimate (\ref{smooth_inc}).

Now we have to analyse the different terms of the above recursion scheme.
For the first term we have
\begin{align}
\vartheta'(\XM_k) (\XM_{k+1} - \XM_k) &=  \frac{a}{\sigma}(\XM_{k}) \label{rec_ym_1}
\Delta + \Delta_k B + \frac{1}{2} \sigma'(\XM_k) (\Delta_k B)^{2} \\ & \nonumber 
\qquad +\frac{1}{2}\left(a'+ \frac{a\sigma'}{\sigma} \right)(\XM_k) \Delta_k B \Delta
+ \frac{1}{2} \frac{aa'}{\sigma}(\XM_k) \Delta^{2} \\ & \qquad + \frac{1}{6} (\sigma \sigma'' + (\sigma')^{2})(\XM_k)(\Delta_k B)^{3}, \nonumber  \end{align} 
while for the second term we obtain
\begin{align} \label{rec_ym_2}
\frac{1}{2}\vartheta''(\XM_k) (\XM_{k+1} - \XM_k)^{2}= - &
\frac{1}{2}\frac{a^{2}\sigma'}{\sigma^{2}}(\XM_k) \Delta^{2}
-\frac{1}{2}\sigma'( \XM_k  ) \left( \Delta_k B \right)^{2} \\ & - \frac{a \sigma'}{\sigma}( \XM_k  )  \Delta_k B \Delta - \frac{1}{2}(\sigma')^{2}( \XM_k  )
(\Delta_k B)^3 + R_k^{(2)} \nonumber
\end{align} 
with
$$   \sup_{k=0, \ldots, n-1}  \EX | R_k^{(2)}|^{p} \leq c \cdot \Delta^{p(1+2H)}. $$
Finally, for the third term we have
\begin{align}\label{rec_ym_3}
\frac{1}{6}\vartheta'''(\XM_k) (\XM_{k+1} - \XM_k)^{3}=
-\frac{1}{6}\left( \sigma'' \sigma - 2(\sigma')^{2}\right)(\XM_k) 
(\Delta_k B)^3 + R_k^{(3)},
\end{align} 
where
 $$   \sup_{k=0, \ldots, n-1} \EX| R_k^{(3)}|^{p} \leq c \cdot \Delta^{p(1+2H)}.$$
Combining  (\ref{rec_ym_1}), (\ref{rec_ym_2}) and   (\ref{rec_ym_3}), it follows
 \begin{align*}
\YM_{k+1}
         = \YM_k+ \Delta_k B + \left( \frac{a}{\sigma} \right) (\XM_k)
         \Delta & +
         \frac{1}{2}\left( a- \frac{a \sigma'}{\sigma} \right) (\XM_k)
         \Delta_k B \Delta  \\ & +  \frac{1}{2}\left( \frac{aa'}{\sigma} -
           \frac{a^{2} \sigma'}{\sigma^{2}} \right)  (\XM_k) \Delta^{2}
         + R_k^{(4)}
\end{align*}        
with $$ \sup_{k=0, \ldots, n-1} \EX| R_k^{(4)}|^{p} \leq c \cdot \Delta^{p(1+2H)}.$$
Using Lemma  \ref{lemma_g} we can write the above recursion as
\begin{align}\label{rec_ym_f}
\YM_{k+1}
         = \YM_k+ \Delta_k B + g(\YM_k)
         \Delta & +
         \frac{1}{2}g' (\YM_k)
         \Delta_k B \Delta  +  \frac{1}{2}gg'  (\YM_k) \Delta^{2}
         + R_k^{(4)}
\end{align}
with $k=0, \ldots, n-1$.

(ii) Now it remains to analyse the difference between $\YL$ and $\YM$.
For this set $t_i= i \Delta$ for $i=0, \ldots, n$. Moreover,
recall that $\YL$ is given by   
$$ \YL_t =\vartheta(\XL_t), \qquad t \in [0,1],$$
and that $\YL$ satisfies the integral  equation
$$ \YL_t =    \YL_{t_i} + \int_{t_i}^{t}  g(\YL_s) \, ds + \Delta_i B \frac{s-t_i}{\Delta} 
 \qquad t \in [t_i,t_{i+1}],$$
with $\YL_0= \vartheta(x_0)$.
A Taylor expansion yields
\begin{align} \YL_{t_{k+1}} = \YL_{t_k} +    \Delta_k B + g(\YL_{t_k})
         \Delta  +
         \frac{1}{2}g' (\YL_{t_k})
         \Delta_k B \Delta +  \frac{1}{2}gg'  (\YL_{t_k}) \Delta^{2}
         + R_k^{(5)}     \label{taylor_yl} \end{align}
for $k=0, \ldots, n-1$ with 
\begin{align*}
 R_k^{(5)} &= \int_{t_k}^{t_{k+1}} \int_{t_k}^{t} \int_{t_k}^{\tau}
 g(g'g)'(\YL_u) \, du \, d \tau \,d t \\ & \quad +  \int_{t_k}^{t_{k+1}} \int_{t_k}^{t} \int_{t_k}^{\tau}
 \frac{\Delta_k B}{\Delta} (g'g)'(\YL_u) \, du \, d \tau \,d t \\
 & \quad +
\int_{t_k}^{t_{k+1}} \int_{t_k}^{t} \int_{t_k}^{\tau}
 \frac{\Delta_k B}{\Delta} g''g(\YL_u) \, du \, d \tau \,d t  \\ &
 \quad +
\int_{t_k}^{t_{k+1}} \int_{t_k}^{t} \int_{t_k}^{\tau}
 \frac{(\Delta_k B)^{2}}{\Delta^{2}} g''(\YL_u) \, du \, d \tau \,d t. 
\end{align*}
Straightforward computations and the boundedness of $g$, $g'$ and $g''
$ yield
$$ \sup_{k=0, \ldots, n-1} \EX| R_k^{(5)}|^{p} \leq c \cdot \Delta^{p(1+2H)}.$$

(3) Now set $e_k= \YL_{t_k}- \YM_{k}$ for $k=0, \ldots, n$.
Using (\ref{rec_ym_f}) and (\ref{taylor_yl})
 we have
\begin{align*}
e_{k+1} = e_k  &+ ( g(\YL_{t_k}) -g(\YM_{k}) )
\Delta  \\ & +
         \frac{1}{2} ( g' (\YL_{t_k}) -  g' (\YM_{k}) )
         \Delta_k B \Delta +  \frac{1}{2} ( gg'  (\YL_{t_k}) -   gg'  (\YM_{k}) ) \Delta^{2}
         + R_k^{(6)}  
\end{align*}
with $e_0=0$ and
\begin{align} \sup_{k=0, \ldots, n-1} \EX| R_k^{(6)}|^{p} \leq c \cdot \Delta^{p(1+2H)}. \label{final_est} \end{align}
Since $g$, $g'$ and $gg'$ are Lipschitz continuous due to Lemma \ref{lemma_g}, it follows that
\begin{align*}
|e_{k+1}|  \leq |e_k|( 1 + c \Delta   + c \sup_{t \in [0,t]} |B_t| \Delta|) +   |R_k^{(6)}| 
\end{align*}
for $k=0, \ldots, n-1$.
The discrete version of Gronwall's Lemma  yields
\begin{align*}
|e_{n}|  \leq  \exp( c (1+  \sup_{t \in [0,t]} |B_t| )) \sum_{k=0}^{n-1} |R_k^{(6)}|
\end{align*}
and (\ref{final_est}) implies that
$$ \EX |e_{n}|^{p} \leq c \cdot \Delta^{2pH}.$$
Now  (\ref{to_show_cons})  follows from  
$$     |\XL_1 - \XM_n| = |\varphi^{-1}(\YL_1)     -  \varphi^{-1}(\YM_n)| \leq c \cdot |e_n|.$$
\fin
 
\appendix 
\section{Appendix} \label{app}
In this section, we show Proposition \ref{int_problem}  following section 3.5.3 in
\cite{diss}, where the case $\sigma=const.$ is studied.

Recall that we have to show 
 $$  \lim_{n \rightarrow \infty} \, n^{2H+1} \EX \left| \int_{0}^{1}  \mathcal{Y}_t (B_t- \BL_t) \, ds \right|^{2} = |\zeta(-2H)|    \int_{0}^{1} \EX \left| \mathcal{Y}_t \right|^{2}  dt . $$
For this, we will need the Malliavin derivative of $
\mathcal{Y}=(\mathcal{Y}_t)_{t \in [0,1]} $ and some technical Lemmata.
Note that the weight function $\mathcal{Y}$ can be written as
\begin{align*}
 \mathcal{Y}_t &= \sigma(X_1) \left( a' - \frac{a \sigma'}{\sigma}
 \right) (X_t) \exp \left(     \int_{t}^{1}  \left( a' - \frac{a \sigma'}{\sigma}
 \right) (X_{\tau})  \,  d \tau \right), \qquad t \in [0,1],
\end{align*}
see  Proposition \ref{rep_mall},
and that
$$ D_sX_t 
= \sigma(X_t) \exp \left(     \int_{s}^{t}  \left( a' - \frac{a \sigma'}{\sigma}
 \right) (X_{\tau})  \,  d \tau \right)1_{[0,t]}(s), \qquad s,t \in [0,1].$$
Moreover, denote 
$$ h(x)= \left(a' - \frac{a \sigma'}{\sigma} \right)(x), \qquad x \in \R
$$
for notational simplicity.

The next Lemma follows by straightforward computations and the assumptions on $a$ and $\sigma$.

\medskip

\begin{lem}\label{prop_weight}
For all $p \geq 1$ there exists a constant $K(p)>0$ such that
$$ \EX |\mathcal{Y}_t - \mathcal{Y}_s|^{p} \leq K(p) |t-s|^{pH}$$
for all
  $s,t  \in [0,1]$.
\end{lem}

\medskip

 The following Lemma 
can be shown by using the product
rule  (\ref{prod_rule}) and
the chain rule (\ref{chain_rule_md}) in Section \ref{pre_1} for the Malliavin derivative.

 \medskip

\begin{lem}  \label{mal_deriv_weight} We have $\mathcal{Y}_t \in \sk^{1,2}$ for all
  $t  \in [0,1]$
with
\begin{align*}
D_{s} \mathcal{Y}_t=   \left(  \frac{\sigma'}{\sigma} \right)(X_1) \mathcal{Y}_t D_s X_1 
 + h'(X_t) D_s X_t       D_sX_1 + \mathcal{Y}_t \int_{t}^{1}
 h'(X_{\tau}) D_s X_{\tau} \, d \tau
\end{align*}
for $s,t \in [0,1]$.
In particular,  there exists a constant $K>0$ such that
$$ \sup_{s,t \in [0,1]} \left| D_{s} \mathcal{Y}_{t}  \right| < K. $$
\end{lem}

\medskip

Denote in the following $\Delta =1/n$ and
\begin{align*}
t_i &= i \Delta, \qquad i=0, \ldots, n, \\
 t_{i+1/2}&= \frac{1}{2} (t_{i}+t_{i+1}), \qquad  i=0, \ldots, n-1.
\end{align*}
We will also use  the notation $\mathcal{Y}_i$ instead of $\mathcal{Y}_{t_i}$
and $D_s \mathcal{Y}_i$ instead of $D_s \mathcal{Y}_{t_i}$ in what follows.
Moreover, set
$$\phi(s,t)=H(2H-1)\cdot  |s-t|^{2H-2}, \qquad s,t \in [0,1].$$
Finally, we will again denote constants, which depend only on  $x_0$, $H$ and $a$, $\sigma$ and their derivatives by $c$, regardless of their value.

\medskip

\begin{lem}{\label{technische_rettung}} We have
$$  \sup_{i,j=0, \ldots, n-1} \, \left|  \int_{t_{j}}^{t_{j+1}} \int_{0}^{1} D_{s}\mathcal{Y}_{i} \cdot (t-t_{j+1/2})  \cdot \phi(s,t) \, ds \, dt \right| \, \leq \,  c \cdot \Delta^{2H+1} .$$
\end{lem}

\medskip

{ \noindent \bf{Proof.
}}  We have
\begin{align} {\label{decomp_tech_lem}}
 &\int_{t_{j}}^{t_{j+1}} \int_{0}^{1} D_{s}\mathcal{Y}_{i}  \cdot (t-t_{j+1/2})  \cdot \phi(s,t) \, ds \, dt  \\ &\qquad  = \int_{t_{j}}^{t_{j+1}} \int_{0}^{t_{j-1}} D_{s}\mathcal{Y}_{i} ¸\cdot   (t-t_{j+1/2})  \cdot \phi(s,t) \, ds \, dt \nonumber \\ &\qquad \qquad  + \int_{t_{j}}^{t_{j+1}} \int_{t_{j-1}}^{t_{j+2}} D_{s}\mathcal{Y}_{i} \cdot  (t-t_{j+1/2})  \cdot \phi(s,t) \, ds \, dt \nonumber
\\  &\qquad \qquad  + \int_{t_{j}}^{t_{j+1}} \int_{t_{j+2}}^{1} D_{s}\mathcal{Y}_{i}  \cdot (t-t_{j+1/2})  \cdot \phi(s,t) \, ds \, dt, \nonumber 
\end{align} with the convention that $t_{-1}=0$ and $t_{n+1}=1$.
 
(i) We start by considering  the second integral. Here we obtain   \begin{align*}
& \left| \int_{t_{j}}^{t_{j+1}} \int_{t_{j-1}}^{t_{j+2}} D_{s}\mathcal{Y}_{i} \cdot (t-t_{j+1/2})  \cdot \phi(s,t) \, ds \, dt \right| 
 \\ & \qquad \qquad \qquad \leq  c \cdot \Delta  \cdot \int_{t_{j}}^{t_{j+1}} \int_{t_{j-1}}^{t_{j+2}}  \phi(s,t) \, ds \, dt
\end{align*}
using Lemma \ref{mal_deriv_weight}.
Since moreover  $$ \int_{t_{j}}^{t_{j+1}} \int_{t_{j-1}}^{t_{j+2}}  \phi(s,t) \, ds \,dt  \leq c \cdot  \Delta^{2H},$$
by straightforward calculations, it follows
\begin{align*}
&\left| \int_{t_{j}}^{t_{j+1}} \int_{t_{j-1}}^{t_{j+2}} D_{s} \mathcal{Y}_{i} \cdot (t-t_{j+1/2})  \cdot \phi(s,t) \, ds \, dt \right| 
 \leq  c \cdot \Delta^{2H+1}. 
\end{align*}

(ii) Now we study the first integral of the right hand side of (\ref{decomp_tech_lem}). This integral clearly vanishes, if $j=0$ or $j=1$. Now  consider $j=2, \ldots, n-1$. In this case
$ 0\leq  s \leq t_{j-1} \leq t_{j} \leq t  \leq t_{j+1}$ and thus we  
have
\begin{align*}
\int_{t_{j}}^{t_{j+1}}  (t-t_{j+1/2}) \cdot \phi(t,s) \, dt   <  0
\end{align*}
for all $s \in [0,t_{j+1}]$. The mean value theorem now implies that
\begin{align*}
& \int_{t_{j}}^{t_{j+1}} \int_{0}^{t_{j-1}} D_{s}\mathcal{Y}_{i} \cdot (t-t_{j+1/2})  \cdot \phi(s,t) \, ds \, dt \\ & \qquad =
  \int_{0}^{t_{j-1}}  D_{s}\mathcal{Y}_{i}  \int_{t_{j}}^{t_{j+1}}  (t-t_{j+1/2})  \cdot \phi(t,s) \, dt \, ds
 \\ & \qquad =
\mu_{i,j} \cdot  \int_{t_{j}}^{t_{j+1}} \int_{0}^{t_{j-1}} (t-t_{j+1/2})  \cdot \phi(s,t) \, ds \, dt ,
\end{align*} where
 $\mu_{i,j}$ is random and satisfies
$$  \inf_{s \in [0,t_{j-1}]} D_{s} \mathcal{Y}_{i}      \,  \leq \,  \mu_{i,j} \, \leq \, \sup_{s \in [0,t_{j-1}]} D_{s} \mathcal{Y}_{i}. $$

By partial integration we obtain
\begin{align*} 
&\int_{t_{j}}^{t_{j+1}}  (t-t_{j+1/2}) \cdot \phi(t,s) \, dt \nonumber \\ & \quad =  \frac{H}{2} \Delta \cdot \left( (t_{j+1}-s)^{2H-1} + (t_{j}-s)^{2H-1} \right) - H  \int_{t_{j}}^{t_{j+1}}  (t-s)^{2H-1} \, dt  \nonumber \\ & \quad = 
\frac{H}{2}  \int_{t_{j}}^{t_{j+1}}  (t_{j+1}-s)^{2H-1} + (t_{j}-s)^{2H-1}  -  2(t-s)^{2H-1} \, d t. \
\end{align*}
Integrating  with respect to the variable $s$ now yields
\begin{align*}
& \int_{t_{j}}^{t_{j+1}} \int_{0}^{t_{j-1}}  (t - t_{j+1/2}) \cdot \phi(s,t)  \, ds \, dt \\& \qquad = -\frac{\Delta}{4} \left( (t_{j+1}-t_{j-1})^{2H} - t_{j+1}^{2H}  +  (t_{j}-t_{j-1})^{2H} - t_{j}^{2H} \right)  \\ & \quad \qquad + \frac{1}{2} \int_{t_{j}}^{t_{j+1}} (t-t_{j-1})^{2H} -t^{2H} \, dt.
\end{align*} 
Note that the above expression is  the error of the trapezoidal approximation for
the integral 
$$ \int_{t_{j}}^{t_{j+1}} g_{j}(t) \, dt $$
with $$g_{j}: [t_{j},t_{j+1}] \rightarrow \R, \qquad g_{j}(t) = \frac{1}{2} \cdot\left(
(t-t_{j-1})^{2H} -t^{2H} \right),  \quad t \in [t_{j},t_{j+1}]. $$
Since $g_{j} \in C^{2}([t_{j},t_{j+1}])$, it clearly holds
\begin{align*}
& \left| \int_{t_{j}}^{t_{j+1}} \int_{0}^{t_{j-1}}  (t - t_{j+1/2}) \cdot \phi(s,t)  \, ds \, dt  \right| =  \frac{g_{j}''(\xi_{j})}{12} \cdot \Delta^{3}
\end{align*}
with $\xi_{j} \in (t_{j},t_{j+1})$. Since
$$g''(t)= 2H(2H-1) \cdot \left( (t-t_{j-1})^{2H-2}-t^{2H-2} \right), \qquad t \in (t_{j},t_{j+1}),$$ 
 we obtain
$$ \sup_{j=2, \ldots, n-1} |  g_{j}''(\xi_{j})|  \leq c \cdot \Delta^{2H-2}.$$

So finally, it follows
$$ \left|  \int_{t_{j}}^{t_{j+1}} \int_{0}^{t_{j-1}} D_{s}\mathcal{Y}_{i} \cdot (t-t_{j+1/2})  \cdot \phi(s,t) \, ds \, dt \right|
     \leq c \cdot \Delta^{1+2H}. $$

(iii) It remains to consider the third integral, which vanishes for $j=n-2$ and $j=n-1$. Here we have $0  \leq  t_{j} \leq t \leq  t_{j+1} \leq t_{j+2} \leq s \leq 1$ and thus
\begin{align*}
\int_{t_{j}}^{t_{j+1}}  (t-t_{j+1/2}) \cdot \phi(s,t) \, dt   > 0 \end{align*}
for all $s \in [t_{j+2},1]$.
Hence we get by the  mean value theorem for integration 
\begin{align*}
& \int_{t_{j}}^{t_{j+1}} \int_{t_{j+2}}^{1} D_{s}\mathcal{Y}_{i} \cdot (t-t_{j+1/2})  \cdot \phi(s,t) \, ds \, dt \\ & \qquad =  \mu_{i,j} \cdot 
\int_{t_{j}}^{t_{j+1}} \int_{t_{j+2}}^{1} (t-t_{j+1/2})  \cdot \phi(s,t) \, ds \, dt ,
\end{align*} with $$  \inf_{s \in [t_{j+2},1]}D_{s}\mathcal{Y}_{i} \, \leq \,  \mu_{i,j} \, \leq \, \sup_{s \in [t_{j+2},1]}D_{s}\mathcal{Y}(t_{i}) .$$
Now, the term $$  \int_{t_{j}}^{t_{j+1}} \int_{t_{j+2}}^{1}  (t-t_{j+1/2})  \cdot \phi(s,t) \, ds  dt     $$ can be treated analogous to (ii).
\fin

Define $$w_{j}(t)= t -t_{j+1/2}, \qquad t \in [t_{j},t_{j+1}],$$
for $j=0, \ldots, n-1$.
In the next Lemma, we precompute the error of the weighted integration problem.

\medskip
\begin{lem}{\label{lem_precomp}}
We have
\begin{align*}
 & \left( \EX \left|\int_{0}^{1} 
(B_t-\BL_t) \mathcal{Y}_t \, dt \right|^{2} \right)^{1/2} \\ & \quad =  \left( \sum_{i=0}^{n-1}\sum_{j=0}^{n-1}  \EX \, \mathcal{Y}_{i}\mathcal{Y}_{j} \int_{t_{i}}^{t_{i+1}}\int_{t_{j}}^{t_{j+1}}  w_{j}(t) w_{i}(s) \phi(s,t) \, dt \, ds \right)^{1/2} +  O(\Delta^{2H}). 
\end{align*}
\end{lem}

\medskip

 {\noindent \bf{Proof.}} First note that 
\begin{align*} & \left( \EX \left|\int_{0}^{1} 
(B_s-\BL_s)\mathcal{Y}_s \, ds \right|^{2} \right)^{1/2} \\ &  \qquad \qquad = \left( \EX \left|\sum_{i=0}^{n-1}\mathcal{Y}_{i} \int_{t_{i}}^{t_{i+1}} 
(B_s-\BL_s) \, ds \right|^{2} \right)^{1/2} + O(\Delta^{2H})
\end{align*} by Lemma \ref{prop_weight}. Moreover, applying partial integration we get
\begin{align*}
 \int_{t_{j}}^{t_{j+1}} 
B_t-\BL_t \, dt & = - \int_{t_{j}}^{t_{j+1}} w_{j}(t) \, d B_t
\end{align*} for $j=0, \ldots, n-1$.
By relation (\ref{rel}) in Subsection  \ref{pre_1} and 
Lemma  \ref{mal_deriv_weight} 
it follows
\begin{align*}
 \mathcal{Y}_{j}\int_{t_{j}}^{t_{j+1}} w_{j}(t) \, d B_{t} & =  
\int_{t_{j}}^{t_{j+1}} \int_{0}^{1} D_{s}\mathcal{Y}_{j} w_{j}(t) \phi(s,t) \, ds \, dt  + \delta ( \mathcal{Y}_j w_j 1_{[t_j,t_{j+1}]}).
\end{align*}
 By Lemma {\ref{technische_rettung}} we have
$$\left| \int_{t_{j}}^{t_{j+1}} \int_{0}^{1} D_{s}\mathcal{Y}_{j} w_{j}(t)\phi(s,t) \, ds \, dt  \right| \leq c \cdot \Delta^{2H+1}.$$ 
Hence, it follows
\begin{align*}
& \left(\EX \left|\sum_{i=0}^{n-1}\mathcal{Y}_{i} \int_{t_{i}}^{t_{i+1}} 
(B_s-\BL_s) \, ds \right|^{2} \right)^{1/2}
\\ & \qquad  \qquad =  \left( \EX \left|\sum_{i=0}^{n-1}  \delta ( \mathcal{Y}_j w_j 1_{[t_j,t_{j+1}]}) \right|^{2} \right)^{1/2} + O(\Delta^{2H}).
\end{align*}
Using the isometry  (\ref{iso_work}) for Skorohod integrals with respect to fractional Brownian motion, we obtain
\begin{align*}
& \EX \left(  \delta ( \mathcal{Y}_i w_i 1_{[t_i,t_{i+1}]})  \delta ( \mathcal{Y}_j w_j 1_{[t_j,t_{j+1}]}) \right)
 \\ & \qquad = 
 \EX \int_{t_{i}}^{t_{i+1}} \int_{t_{j}}^{t_{j+1}} \mathcal{Y}_{i} \mathcal{Y}_{j} \cdot  w_{i}(s) w_{j}(t)  \cdot \phi(s,t)  \, dt \, ds \\ &  \qquad \qquad \quad + 
\EX \int_{t_{i}}^{t_{i+1}} \int_{t_{j}}^{t_{j+1}} \int_{0}^{1} \int_{0}^{1} D_{s_{1}}\mathcal{Y}_{i} w_{i}(t_{1})  \cdot D_{s_{2}} \mathcal{Y}_{j}w_{j}(t_{2})  \\ & \qquad \qquad \qquad \qquad \qquad \qquad \qquad \qquad \qquad \times \phi(s_{1},t_{2}) \phi(s_{2},t_{1}) \, d s_{2} \, ds_{1} \, d t_{2} \, d t_{1}
\end{align*} for $i,j=0, \ldots, n-1$.
Since fortunately
\begin{align*}
& \int_{t_{i}}^{t_{i+1}} \int_{t_{j}}^{t_{j+1}} \int_{0}^{1} \int_{0}^{1} D_{s_{1}} \mathcal{Y}_{i} w_{i}(t_{1})  \cdot D_{s_{2}} \mathcal{Y}_{j} w_{j}(t_{2})  \cdot \phi(s_{1},t_{2}) \phi(s_{2},t_{1}) \, d s_{2} \, ds_{1} \, d t_{2} \, d t_{1} \\ &=  \int_{t_{i}}^{t_{i+1}} \int_{0}^{1} D_{s_{2}}\mathcal{Y}_{j} w_{i}(t_{1})  \phi(s_{2},t_{1}) \, ds_{2} \, dt_{1}    \cdot   \int_{t_{j}}^{t_{j+1}} \int_{0}^{1} D_{s_{1}}\mathcal{Y}_{i} w_{j}(t_{2}) \phi(s_{1},t_{2}) \, ds_{1} \, dt_{2}, 
\end{align*}
we  have by Lemma {\ref{technische_rettung}} that
\begin{align*}
& \EX \left(  \delta ( \mathcal{Y}_i w_i 1_{[t_i,t_{i+1}]})  \delta ( \mathcal{Y}_j w_j 1_{[t_j,t_{j+1}]}) \right)
\\ & \qquad = 
 \EX \int_{t_{i}}^{t_{i+1}} \int_{t_{j}}^{t_{j+1}} \mathcal{Y}_{i} \mathcal{Y}_{j} \cdot  w_{i}(s) w_{j}(t)  \cdot \phi(s,t)  \, dt \, ds  + O(\Delta^{4H+2}),
\end{align*}  
which shows the assertion. \fin

Now we finally determine the strong asymptotic behaviour  of
$$\EX \left| \int_{0}^{1}  \mathcal{Y}_t (B_t- \BL_t) \, dt \right|^{2}.$$
 For similar calculations in the case that the weight function $\mathcal{Y}$ is deterministic  and  the  process $B$ is stationary and  behaves locally like a fractional Brownian motion, see e.g.  \cite{stein} and \cite{benh}.

{ \noindent \bf Proof of  Proposition  \ref{int_problem}.}
By Lemma 
{\ref{lem_precomp}}
we have
\begin{align*}
 & \left( \EX \left|\int_{0}^{1}  \mathcal{Y}_t (B_t- \BL_t) \, dt \right|^{2} \right)^{1/2} \\ & \quad =  \left( \sum_{i=0}^{n-1}\sum_{j=0}^{n-1}  \EX \, \mathcal{Y}_{i}\mathcal{Y}_{j} \int_{t_{i}}^{t_{i+1}}\int_{t_{j}}^{t_{j+1}}  w_{i}(s) w_{j}(t) \cdot  \phi(s,t) \, dt \, ds \right)^{1/2} +  O(\Delta^{2H}). 
\end{align*}
Thus, it remains to study the behaviour of
$$  \sum_{i=0}^{n-1} \sum_{j=0}^{n-1} \int_{t_{i}}^{t_{i+1}} \int_{t_{j}}^{t_{j+1}} \EX  \, \mathcal{Y}_{i} \mathcal{Y}_{j} \cdot  w_{i}(s) w_{j}(t)  \cdot \phi(s,t)  \, dt \, ds.$$

Note that 
\begin{align*}
  \int_{t_{i}}^{t_{i+1}}  \int_{t_{j}}^{t_{j+1}} w_{i}(s)w_{j}(t) \cdot \phi(s,t) \, dt \, ds   = \EX \int_{t_{i}}^{t_{i+1}}   w_{i}(s) \, dB_s    \int_{t_{j}}^{t_{j+1}}  w_{j}(t) \, dB_t
\end{align*}  by  (\ref{iso_work}) in Section \ref{pre_1}, and recall that
\begin{align*}
& \int_{t_{i}}^{t_{i+1}}  w_{i}(s) \, dB_s   = \int_{t_{i}}^{t_{i+1}} \frac{1}{2} \left( B_{t_{i}} +B_{t_{i+1}} \right)-B_s \, ds.
\end{align*} 
Define 
\begin{align*}
 & \theta_{i,j}(s_{1},s_{2})   =\frac{1}{4}\EX ( B_{t_{i}} +B_{t_{i+1}} -2B_{s_{1}})( B_{t_{j}} +B_{t_{j+1}} -2B_{s_{2}})
\end{align*}
for $s_{1} \in [t_{i},t_{i+1}]$, $s_{2} \in [t_{j},t_{j+1}]$, $i,j = 0, \ldots, n-1.$
Thus we can write
\begin{align*}  & \sum_{i=0}^{n-1} \sum_{j=0}^{n-1} \int_{t_{i}}^{t_{i+1}} \int_{t_{j}}^{t_{j+1}} \EX  \, \mathcal{Y}_{i} \mathcal{Y}_{j} \cdot  w_{i}(s) w_{j}(t)  \cdot \phi(s,t)  \, dt \, ds \\ & \qquad \qquad  = 
 \sum_{i=0}^{n-1} \sum_{j=0}^{n-1} \int_{t_{i}}^{t_{i+1}} \int_{t_{j}}^{t_{j+1}} \EX  \, \mathcal{Y}_{i} \mathcal{Y}_{j} \cdot  \theta_{i,j}(s_{1},s_{2})  \, ds_{2} \, ds_{1}. 
\end{align*}
By straightforward calculations we obtain
\begin{align*}
  \theta_{i,j}(s_{1},s_{2}) &= -\frac{1}{8} \left( |t_{i}-t_{j}|^{2H} + |t_{i}-t_{j+1}|^{2H} + |t_{i+1}-t_{j}|^{2H} + |t_{i+1}-t_{j+1}|^{2H} \right)\\ &   \quad + 
 \frac{1}{4}\left( |t_{i}-s_{2}|^{2H} + |t_{i+1}-s_{2}|^{2H} + |t_{j}-s_{1}|^{2H} + |t_{j+1}-s_{1}|^{2H} \right) \\ &   \quad-\frac{1}{2} |s_{1}-s_{2}|^{2H},
\end{align*}
which simplifies in the case $i=j$ to
\begin{align*}
   \theta_{j,j}(s_{1},s_{2}) & = -\frac{1}{4}  |t_{j+1}-t_{j}|^{2H} \\ &   \quad + 
 \frac{1}{4}\left( |t_{j}-s_{2}|^{2H} + |t_{j+1}-s_{2}|^{2H} + |t_{j}-s_{1}|^{2H} + |t_{j+1}-s_{1}|^{2H} \right) \\ &   \quad -\frac{1}{2} |s_{1}-s_{2}|^{2H}.
\end{align*}

(i) We first show that asymptotically the contribution of the off-diagonal terms  to the error is negligible, i.e., 
\begin{align}\label{off-diag}
 \sum_{|i-j| > \log(n)}\int_{t_{i}}^{t_{i+1}} \int_{t_{j}}^{t_{j+1}} \EX \mathcal{Y}_{i}\mathcal{Y}_{j} \cdot \theta_{i,j}(s_{1},s_{2}) \, ds_{2} \, ds_{1} = o(n^{-2H-1}).\end{align}
Note that by symmetry
\begin{align*}
& \sum_{|i-j| > \log(n)}\int_{t_{i}}^{t_{i+1}} \int_{t_{j}}^{t_{j+1}} \EX  \mathcal{Y}_{i}\mathcal{Y}_{j} \cdot \theta_{i,j}(s_{1},s_{2}) \, ds_{2} \, ds_{1} \\ & \qquad  =2  \cdot  \sum_{i-j > \log(n)}\int_{t_{i}}^{t_{i+1}} \int_{t_{j}}^{t_{j+1}} \EX \mathcal{Y}_{i}\mathcal{Y}_{j} \cdot \theta_{i,j}(s_{1},s_{2}) \, ds_{2} \, ds_{1}
\end{align*}
To show (\ref{off-diag}) we will use  fourth order Taylor expansions  
of suitable parts of $\theta_{i,j}(s_{1},s_{2})$. 
For this, the following  will be very helpful:

Let $a>0$ and $\epsilon_{x}, \epsilon_{y} \in \{ -1, +1 \}$. Then  for $x,y \in [0,1]$ such that $a+ \epsilon_{x} x +  \epsilon_{y}y >0$, the function
$$ f(x,y)=(a+\epsilon_{x} x +  \epsilon_{y}y)^{2H}$$ is well defined and we have \begin{align} {\label{helpful_form}}\frac{\partial^{n}f}{(\partial x)^{k} (\partial y)^{n-k}} (x,y) = \kappa_{n} \cdot \epsilon_{x}^{k} \epsilon_{y}^{n-k} \cdot (a+\epsilon_{x} x +  \epsilon_{y}y)^{2H-n} \end{align}
with $$\kappa_{n} = 2H \cdot (2H-1) \cdot \ldots \cdot (2H-n+1).$$

In what follows, set $$ \tau_{i,j}=|t_{i+1/2}-t_{j+1/2}|=\Delta |i-j|,$$ and  recall that $$w_{i}(s_{1})=s_{1}-t_{i+1/2}, \qquad w_{j}(s_{2})=s_{2}-t_{j+1/2}.$$

The first part of $\theta_{i,j}$ we study is 
\begin{align*} \theta_{i,j}^{(1)} & =-\frac{1}{8} \left( |t_{i}-t_{j}|^{2H} + |t_{i}-t_{j+1}|^{2H} + |t_{i+1}-t_{j}|^{2H} + |t_{i+1}-t_{j+1}|^{2H} \right) \\
& =-\frac{1}{8} \left( 2 \tau_{i,j}^{2H} + |\tau_{i,j}-\Delta/2 - \Delta/2|^{2H} \ + |\tau_{i,j}+\Delta/2 + \Delta/2|^{2H}  \right). 
\end{align*}
Since $i>j$, we obtain by applying  (\ref{helpful_form}) with $a=\tau_{i,j}$, $x =\Delta/2$ and $y=\Delta/2$  the expansion
$$  \theta_{i,j}^{(1)}= -\frac{1}{2} \cdot \tau_{i,j}^{2H}  -\frac{\kappa_{2}}{8} \cdot \tau_{i,j}^{2H-2} \cdot \Delta^{2}  + \rho_{i,j},$$
with
$$ |\rho_{i,j}| \leq  c \cdot |t_{i}-t_{j+1}|^{2H-4} \cdot \Delta^{4}.$$
Since   $i-j-1 \geq \log(n)$, we have
$$  \rho_{i,j} = O\left( \log(n)^{4-2H}n^{-4} \right).$$

Now consider the second part of $\theta_{i,j}$ given by
\begin{align*}
\theta_{i,j}^{(2)}(s_{1},s_{2}) & =\frac{1}{4} \left((t_{i}-s_{2})^{2H} +(t_{i+1}-s_{2})^{2H} \right) \\ & =\frac{1}{4}\left((\tau_{i,j}-  w_{j}(s_{2}) - \Delta/2)^{2H} +(\tau_{i,j}-  w_{j}(s_{2})+ \Delta/2)^{2H} \right). \end{align*}
Here we obtain by applying (\ref{helpful_form}) with $a=\tau_{i,j}$, $x =  w_{j}(s_{2})$ and $y=\Delta/2$ 
\begin{align*}
\theta_{i,j}^{(2)}(s_{1},s_{2}) & =  \frac{1}{2} \cdot \tau_{i,j}^{2H}  + \frac{\kappa_{2}}{4} \cdot \tau_{i,j}^{2H-2} \cdot \left((s_{2}-t_{j+1/2})^{2} + \Delta^{2}/4 \right) \\ & \qquad -  \frac{\kappa_{1}}{2} \cdot  \tau_{i,j}^{2H-1} \cdot  (s_{2}-t_{j+1/2}) \\ & \qquad  - \frac{\kappa_{3}}{12} \cdot   \tau_{i,j}^{2H-3} \cdot (s_{2}-t_{j+1/2} )^{3} \\ &\qquad - \frac{\kappa_{3}}{16}  \cdot  \tau_{i,j}^{2H-3} \cdot  (s_{2}-t_{j+1/2} ) \cdot \Delta^{2}
\\ & \qquad + O\left( \log(n)^{4-2H}n^{-4} \right). 
\end{align*} 
Note that clearly 
$$ \int_{t_{j}}^{t_{j+1}} s_{2}-t_{j+1/2} \, ds_{2} =0, \qquad  \int_{t_{j}}^{t_{j+1}} (s_{2}-t_{j+1/2})^{3} \, ds_{2} =0,$$
so the third, fourth and fifth term of the above expansion vanish after integration over $[t_{i},t_{i+1}] \times [t_{j}, t_{j+1}]$. 
For
\begin{align*}
\theta_{i,j}^{(3)}(s_{1},s_{2}) & =\frac{1}{4} \left((s_{1}-t_{j})^{2H} +(s_{1}-t_{j+1})^{2H} \right) \\
& = \frac{1}{4} \left((\tau_{i,j}+ w_{i}(s_{1}) +\Delta/2)^{2H} 
+(\tau_{i,j}+  w_{i}(s_{1}) -\Delta/2)^{2H}   \right) 
\end{align*}
we obtain analogously 
\begin{align*}
\theta_{i,j}^{(3)}(s_{1},s_{2}) & =  \frac{1}{2}  \cdot \tau_{i,j}^{2H} +  \frac{\kappa_{2}}{4} \cdot  \tau_{i,j}^{2H-2} \cdot \left((s_{1}-t_{i+1/2})^{2} + \Delta^{2}/4 \right)
+ \nu_{i,j} \\ & \qquad +  O\left( \log(n)^{4-2H}n^{-4} \right),
\end{align*}
where $\nu_{i,j}$ denotes the terms of the Taylor expansion, which contain odd powers of $s_{1}-t_{i+1/2}$ resp. $s_{2}-t_{j+1/2}$ and vanish after integration over $[t_{i},t_{i+1}] \times [t_{j}, t_{j+1}]$.
 
Finally, we have to study
$$ \theta^{(4)}_{i,j}(s_{1},s_{2})= -\frac{1}{2} |s_{1}-s_{2}|^{2H}= -\frac{1}{2}|\tau_{i,j} + w_{i}(s_{1}) - w_{j}(s_{2})|^{2H}.$$
Here we obtain
\begin{align*}
 \theta^{(4)}_{i,j}(s_{1},s_{2})& = -\frac{1}{2}\cdot \tau_{i,j}^{2H}  - \frac{\kappa_{2}}{4} \cdot \tau_{i,j}^{2H-2} \cdot \left((s_{1}-t_{i+1/2})^{2}+ (s_{2}-t_{j+1/2})^{2}\right)+\nu_{i,j}  \\ & \qquad + O\left( \log(n)^{4-2H}n^{-4} \right), 
\end{align*}
where $\nu_{i,j}$ denotes again the terms, which vanish after integration.

Summing up the above expansions for the parts of  $\theta_{i,j}$ yields
$$\theta_{i,j}(s_{1},s_{2})=\nu_{i,j} + O\left( \log(n)^{4-2H}n^{-4} \right).$$ Therefore we obtain after integrating over $[t_{i},t_{i+1}] \times [t_{j}, t_{j+1}]$ the estimate
$$ \left| \int_{t_{i}}^{t_{i+1}} \int_{t_{j}}^{t_{j+1}} \theta_{i,j}(s_{1},s_{2}) \, ds_{2} \, ds_{1} \right| \leq c \cdot (\log(n))^{4-2H} \cdot \frac{1}{n^{6}}.$$  Since $\mathcal{Y}$ is bounded, it finally  follows 
\begin{align*}
 \left|  \sum_{|i-j| > \log(n)}\int_{t_{i}}^{t_{i+1}} \int_{t_{j}}^{t_{j+1}} \EX \mathcal{Y}_{i} \mathcal{Y}_{j} \cdot \theta_{i,j}(s_{1},s_{2}) \, ds_{2} \, ds_{1} \right| \leq c \cdot (\log(n))^{4-2H} \cdot \frac{1}{n^{4}}, 
\end{align*} and hence we have shown (\ref{off-diag}).

(ii) Now it remains to consider  the summands with $|i-j| \leq \log(n)$, i.e. the diagonal resp. near diagonal parts. For this, we need to compute the integrals over 
$[t_{i},t_{i+1}] \times [t_{j},t_{j+1}]$ of the four parts of $\theta_{i,j}$. 
Note that by symmetry
\begin{align*}
 & \sum_{|i-j| \leq  \log(n)}\int_{t_{i}}^{t_{i+1}} \int_{t_{j}}^{t_{j+1}} \EX \mathcal{Y}_{i} \mathcal{Y}_{j} \cdot \theta_{i,j}(s_{1},s_{2}) \, ds_{2} \, ds_{1} \\ & \qquad = 2 \cdot \sum_{0 < i-j  \leq \log(n)} \int_{t_{i}}^{t_{i+1}} \int_{t_{j}}^{t_{j+1}} \EX \mathcal{Y}_{i} \mathcal{Y}_{j} \cdot \theta_{i,j}(s_{1},s_{2}) \, ds_{2} \, ds_{1} \\ & \qquad \qquad \, \, \, \, + \sum_{ j=0 }^{n-1} \int_{t_{j}}^{t_{j+1}} \int_{t_{j}}^{t_{j+1}} \EX  \mathcal{Y}_{j}^{2} \cdot \theta_{j,j}(s_{1},s_{2}) \, ds_{2} \, ds_{1}
.
\end{align*}
For $i>j$ we obtain
\begin{align}
 & \int_{t_{i}}^{t_{i+1}} \int_{t_{j}}^{t_{j+1}}  \theta_{i,j}^{(1)}(s_{1},s_{2}) \, ds_{2} \, ds_{1} \nonumber \\ & \quad = -\frac{\Delta^2}{8} \left( |t_{i}-t_{j}|^{2H} + |t_{i}-t_{j+1}|^{2H} + |t_{i+1}-t_{j}|^{2H} + |t_{i+1}-t_{j+1}|^{2H} \right) \nonumber  \\  & \quad = 
-\frac{1}{8}   \cdot \left( 2|i-j|^{2H} + |i+1-j|^{2H} +  |i-1-j|^{2H} \right) \cdot  \frac{1}{n^{2H+2}} \label{int_1}.
\end{align}
Moreover, we have 
\begin{align}
&  \int_{t_{i}}^{t_{i+1}} \int_{t_{j}}^{t_{j+1}}  \theta_{i,j}^{(2)}(s_{1},s_{2}) \, ds_{2}  \, ds_1 \nonumber \\ & \quad = 
\frac{\Delta}{4(2H+1)} \left(|t_{i}-t_{j}|^{2H+1} +|t_{i+1}-t_{j}|^{2H+1}\right)  \nonumber \\ & \qquad \qquad \qquad \qquad - \frac{\Delta}{4(2H+1)} \left( |t_{i}-t_{j+1}|^{2H+1} +|t_{i+1}-t_{j+1}|^{2H+1}\right) \nonumber  \\ 
& \quad =\frac{1}{4(2H+1)} \left( |i+1-j|^{2H+1} -|i-1-j|^{2H+1}  \right)  \cdot  \frac{1}{n^{2H+2}}
\label{int_2}, 
 \\ & \int_{t_{i}}^{t_{i+1}} \int_{t_{j}}^{t_{j+1}}  \theta_{i,j}^{(3)}(s_{1},s_{2}) \, ds_{2} \, ds_{1}  \nonumber \\ & \quad = 
\frac{\Delta}{4(2H+1)} \left(|t_{i+1}-t_{j}|^{2H+1} +|t_{i+1}-t_{j+1}|^{2H+1}\right)  \nonumber 
\\ & \qquad \qquad \qquad \qquad - \frac{\Delta}{4(2H+1)} \left( |t_{i}-t_{j}|^{2H+1} +|t_{i}-t_{j+1}|^{2H+1}\right) \nonumber  
\\ & \quad =\frac{1}{4(2H+1)} \left( |i+1-j|^{2H+1} -|i-1-j|^{2H+1}  \right)  \cdot  \frac{1}{n^{2H+2}}
\label{int_3} 
\end{align}
and
\begin{align} \nonumber
 & \int_{t_{i}}^{t_{i+1}} \int_{t_{j}}^{t_{j+1}}  \theta_{i,j}^{(4)}(s_{1},s_{2}) \, ds_{2} \, ds_{1} \\ & \quad = 
\frac{1}{2(2H+1)(2H+2)} \left(|t_{i+1}-t_{j+1}|^{2H+2} +|t_{i}-t_{j}|^{2H+2}\right)  \nonumber \\ & \qquad \qquad \qquad \qquad - \frac{1}{2(2H+1)(2H+2)} \left(|t_{i}-t_{j+1}|^{2H+2} +|t_{i+1}-t_{j}|^{2H+2}\right) \nonumber
\\ & \quad = \frac{1}{2(2H+1)(2H+2)} \left(2|i-j|^{2H+2} \right. \nonumber  \\ & \qquad \qquad \qquad \qquad \qquad  \left. - |i-1-j|^{2H+2} -|i+1-j|^{2H+2} \right) \cdot \frac{1}{n^{2H+2}} .  \label{int_4}
\end{align}
Moreover, for $i=j$ we have
\begin{align} \nonumber & 
\int_{t_{j}}^{t_{j+1}} \int_{t_{j}}^{t_{j+1}} \theta_{j,j}(s_{1},s_{2}) \, ds_{2} \, ds_{1}\\ & \qquad \qquad =  \left(\frac{1}{2H+1}-\frac{1}{(2H+2)(2H+1)} - \frac{1}{4} \right)\cdot
\frac{1}{n^{2H+2}}.\label{int_i=j}
\end{align}
Hence, combining (\ref{int_1}),  (\ref{int_2}),  (\ref{int_3}) and  (\ref{int_4}) we obtain for $i>j$ that
\begin{align*}
& \int_{t_{i}}^{t_{i+1}} \int_{t_{j}}^{t_{j+1}}  \theta_{i,j}(s_{1},s_{2}) \, ds_{2} \, ds_{1}=  \frac{1}{n^{2H+2}}\cdot  \mathcal{K}_{1}(i,j) 
\end{align*}
 with
\begin{align*} 
\mathcal{K}_{1}(i,j) & = -\frac{1}{8}  \cdot \left( 2|i-j|^{2H} + |i+1-j|^{2H} +  |i-1-j|^{2H} \right) \\ & \qquad  + \frac{1}{2(2H+1)} \cdot \left( |i+1-j|^{2H+1} - |i-1-j|^{2H+1}  \right)
\\ & \qquad +  \frac{1}{2(2H+1)(2H+2)} \cdot  \left( 2|i-j|^{2H+2} -|i-j+1|^{2H+2} -|i-j-1|^{2H+2} \right).
\end{align*}
In the case $i=j$ we have by (\ref{int_i=j}) that
\begin{align*} 
& \int_{t_{i}}^{t_{i+1}} \int_{t_{j}}^{t_{j+1}}  \theta_{j,j}(s_{1},s_{2}) \, ds_{2} \, ds_{1} = 
  \mathcal{K}_{2} \cdot \frac{1}{n^{2H+2}}
\end{align*}
with
$$ \mathcal{K}_{2}= \frac{1}{2H+1}-\frac{1}{(2H+2)(2H+1)} - \frac{1}{4}.$$
Hence it follows
\begin{align*}
 & \sum_{|i-j| \leq \log(n)}\int_{t_{i}}^{t_{i+1}} \int_{t_{j}}^{t_{j+1}} \EX \mathcal{Y}_{i} \mathcal{Y}_{j} \cdot \theta_{i,j}(s_{1},s_{2}) \, ds_{2} \, ds_{1} \\ & \qquad = 2 \cdot \sum_{0 < i-j  \leq \log(n)} \EX \mathcal{Y}_{i} \mathcal{Y}_{j} \cdot \mathcal{K}_{1}(i,j)  \cdot \frac{1}{n^{2H+2}}  +  \sum_{ j=0 }^{n-1} \EX \mathcal{Y}_{j}^{2}  \cdot \mathcal{K}_{2} \cdot \frac{1}{n^{2H+2}}.
\end{align*}

Since by Lemma \ref{prop_weight} it holds
$$ \EX|\mathcal{Y}_{t_{k}}- \mathcal{Y}_{t_{l}}|^{2} \leq c \cdot \log(n)^{2} \cdot \frac{1}{n^{2H}}$$ 
for $|k-l| \leq \log(n)$,
we have
\begin{align} \label{almost_final}
 & \sum_{|i-j| \leq \log(n)}\int_{t_{i}}^{t_{i+1}} \int_{t_{j}}^{t_{j+1}} \EX \mathcal{Y}_{i} \mathcal{Y}_{j} \cdot \theta_{i,j}(s_{1},s_{2}) \, ds_{2} \, ds_{1} \\ \nonumber & \qquad = 2 \cdot \sum_{0 < i-j  \leq \log(n)} \EX  \mathcal{Y}_{j}^{2} \cdot \mathcal{K}_{1}(i,j)  \cdot \frac{1}{n^{2H+2}}  +  \sum_{ j=0 }^{n-1} \EX \mathcal{Y}_{j}^{2}  \cdot \mathcal{K}_{2} \cdot \frac{1}{n^{2H+2}} + o\left( n^{-2H-1} \right). \nonumber 
\end{align}

Rearranging the  terms on the right hand side of the above equation yields
\begin{align*}
& 2 \cdot \sum_{0 < i-j  \leq \log(n)} \EX  \mathcal{Y}_{j}^{2} \cdot \mathcal{K}_{1}(i,j)  \cdot \frac{1}{n^{2H+2}}  +  \sum_{ j=0 }^{n-1} \EX \mathcal{Y}_{j}^{2}  \cdot \mathcal{K}_{2} \cdot \frac{1}{n^{2H+2}} \\
& \qquad =   \sum_{j=0}^{n-1} \EX \mathcal{Y}_{j}^{2} \frac{1}{n^{2H+2}} \cdot \mathcal{C}_{j}(\lfloor \log(n) \rfloor),
\end{align*}
where 
\begin{align*}
\mathcal{C}_{j}(r)& = \mathcal{K}_{2}+2 \sum_{i=j+1}^{j+r} \mathcal{K}_{1}(i,j)
\end{align*}
for $j=0, \ldots, n-1$, $r \in \N$.
Note that
\begin{align*}
2 \sum_{i=j+1}^{j+r} \mathcal{K}_{1}(i,j)&= -\frac{1}{4}  \sum_{k= 1}^{r}  2 k^{2H} + (k+1)^{2H} + (k-1)^{2H} 
\\ & \quad +  \frac{1}{2H+1}\sum_{k=1}^{r} (k+1)^{2H+1} - (k-1)^{2H+1}  
\\& \quad +   \frac{1}{(2H+1)(2H+2)}\sum_{k= 1}^{r} 2k^{2H+2} -(k+1)^{2H+2}-(k-1)^{2H+2}\\
 & = -\frac{1}{4}(r+1)^{2H}+ \frac{1}{4}r^{2H}+  \frac{1}{4}-\sum_{k=1}^{r}  k^{2H}
\\ & \quad +  \frac{1}{2H+1} \left( (r+1)^{2H+1} +r^{2H+1} - 1 \right)
\\& \quad +   \frac{1}{(2H+1)(2H+2)} \left(r^{2H+2}-(r+1)^{2H+2}+1\right)
\end{align*} for $r \in \N$.
Thus we have
$$ \mathcal{C}_{j}(r)=\mathcal{C}_{0}(r)$$
for  $j =0, \ldots, n-1 $ and
\begin{align*} \mathcal{C}_{0}(r)& = -\frac{1}{4}(r+1)^{2H}+ \frac{1}{4}r^{2H}-\sum_{k=1}^{r}  k^{2H}
\\ & \quad +  \frac{1}{2H+1}  \left((r+1)^{2H+1} +r^{2H+1} \right)
\\& \quad +    \frac{1}{(2H+1)(2H+2)} \left(r^{2H+2} -(r+1)^{2H+2}\right).
\end{align*} 
In the next step  we show that 
\begin{align} \label{zeta} \lim_{r \rightarrow \infty}  \mathcal{C}_{0}(r)  = -\zeta(-2H).\end{align}  For this, we again apply Taylor expansions of suitable parts of  $\mathcal{C}_{0}(r)$. We have
\begin{align*}
\frac{1}{(2H+1)(2H+2)}(r+1)^{2H+2} &=  \frac{1}{(2H+1)(2H+2)}r^{2H+2} + \frac{1}{2H+1}r^{2H+1} \\ & \quad +   \frac{1}{2}r^{2H}  + \frac{H}{3}r^{2H-1} + o(1), \\
\frac{1}{2H+1}(r+1)^{2H+1} &=   \frac{1}{2H+1}r^{2H+1} +   r^{2H}  + Hr^{2H-1} + o(1),
\\ \frac{1}{4}(r+1)^{2H}& =\frac{1}{4}r^{2H} + \frac{H}{2}r^{2H-1} + o(1).
\end{align*}
Hence it follows
 \begin{align*} \mathcal{C}_{0}(r)& = -\sum_{k=1}^{r}  k^{2H} + \frac{1}{2H+1}r^{2H+1} + \frac{r}{2}^{2H} + \frac{H}{6}r^{2H-1} + o(1).
\end{align*}
Since 
$$ \sum_{k=1}^{r}  k^{2H} = \zeta(-2H)+ \frac{1}{2H+1}r^{2H+1} + \frac{r}{2}^{2H} + \frac{H}{6}r^{2H-1} + o(1), $$ see e.g. \cite{abram}, we obtain
(\ref{zeta}).  

This yields
\begin{align*}
& \lim_{n \rightarrow \infty} \, n^{1+2H} \sum_{|i-j| \leq  \log(n)}\int_{t_{i}}^{t_{i+1}}  \int_{t_{j}}^{t_{j+1}} \EX \mathcal{Y}_{t_{i}}\mathcal{Y}_{t_{j}} \cdot \theta_{i,j}(s_{1},s_{2}) \, ds_{2} \, ds_{1} \\ & \quad =  \lim_{n \rightarrow \infty}
\mathcal{C}_{0}( \lfloor \log(n) \rfloor ) \cdot \sum_{j=0}^{n-1}  \EX \mathcal{Y}_{t_{j}}^{2} \cdot (t_{j+1}-t_{j})\\ & \quad =  |\zeta(-2H)| \cdot \lim_{n \rightarrow \infty}
 \sum_{j=0}^{n-1}  \EX \mathcal{Y}_{t_{j}}^{2} \cdot (t_{j+1}-t_{j}) \\ & \quad =  |\zeta(-2H)| \cdot \int_{0}^{1} \EX \mathcal{Y}_{t}^{2} \, dt,
\end{align*}
which together with (\ref{off-diag})  shows finally Proposition \ref{int_problem}.
\fin

\end{document}